\documentclass[a4paper]{amsart}

\usepackage[mathscr]{eucal}
\usepackage{epsfig}
\usepackage{amssymb}
\usepackage{longtable}
\usepackage{amssymb}
\usepackage{amscd,amsfonts}
\usepackage{amsthm}
\usepackage{amsmath}
\usepackage{latexsym}

\theoremstyle{plain}
\newtheorem{theorem}{Theorem}

\newtheorem{corollary}[theorem]{Corollary}

\newtheorem{lemma}[theorem]{Lemma}

\theoremstyle{definition}
\newtheorem{definition}[theorem]{Definition}
\newtheorem{remark}[theorem]{Remark}

\newtheorem{notation}[theorem]{Notation}

\numberwithin{theorem}{section}

\newcommand{\Ch}{\mathbb{C}}

\newcommand{\card}{\mbox{card}}

\newcommand{\Proof}{{\it Proof. \,}}
\newcommand{\dgFv}{{\mathit dg}\widetilde{{\mathcal F}_v}}
\newcommand{\barFv}{{\widetilde{{\mathcal F}_v}}}

\newcommand{\barSo}{{\widetilde{{\mathcal S_v}^{\perp}}}}
\newcommand{\Hom}{\mbox{Hom}}
\newcommand{\Uo}{{\mathcal L}^{\perp}}
\newcommand{\Co}{{\mathcal C}^{\perp}}
\newcommand{\oUo}{{^{\perp}}({\mathcal L}^{\perp})}
\newcommand{\R}{{\mathfrak R}}

\def\Qco{\mathfrak{Qco}}
\def\M{\mathscr{M}}
\def\A{\mathscr{A}}

\newcommand{\dgC}{{\mathit dg\,}\widetilde{{\mathcal C}}}
\newcommand{\barC}{{\widetilde{{\mathcal C}}}}
\newcommand{\dgCo}{{\mathit dg\,}\widetilde{{\mathcal C}^{\perp}}}
\newcommand{\barCo}{{\widetilde{{\mathcal C}^{\perp}}}}
\newcommand{\barLo}{{\widetilde{{\mathcal L}^{\perp}}}}

\newcommand{\dgCoo}{{\mathit dg\,}\widetilde{{\mathcal L}^{\perp}}}
\newcommand{\barCoo}{{\widetilde{{\mathcal L}^{\perp}}}}

\def\Qco{\mathfrak{Qco}}
\def\Mod{\mbox{-Mod}}

\def\Z{\mathbb{Z}}

\newcommand{\OHom}[3]{\mbox{Hom}_{#1}(#2,#3)}
\newcommand{\Ext}[4]{\mbox{Ext}^{#1}_{#2}(#3,#4)}

\newcommand{\rmod}[1]{\mbox{\rm{Mod}--}{#1}}

\newcommand{\Spec}{\mbox{\rm Spec}}

\begin{document}

\title{Model category structures arising from Drinfeld vector bundles}
\author{Sergio Estrada}
\address{Departamento de Matem\'atica Aplicada, Universidad de Murcia\\
Campus del Espinardo 30100, Spain}
\author{Pedro Guil Asensio}
\address{Departamento de Matem\'aticas, Universidad de Murcia\\
Campus del Espinardo 30100, Spain}
\author{Mike Prest}
\address{School of Mathematics, University of Manchester\\
Oxford Road, Manchester M13 9PL,  UK}
\author{Jan Trlifaj}
\address{Charles University, Faculty of Mathematics and Physics, Department of Algebra \\
Sokolovsk\'{a} 83, 186 75 Prague 8, Czech Republic}
\thanks{The first two authors have been partially supported by the DGI and by the
Fundaci\'on Seneca.
The fourth author supported by GA\v CR 201/09/0816 and MSM 0021620839}


\date{\today}
\begin{abstract}
We present a general construction of model category structures on the category
$\Ch(\Qco(X))$ of unbounded chain complexes
of quasi-coherent sheaves on a semi-separated scheme $X$. The construction is based
on making compatible the filtrations of individual
modules of sections at open affine subsets of $X$. It does not require closure under
direct limits as previous methods. We apply it to describe the
derived category $\mathbb D (\Qco(X))$
via various model structures on $\Ch(\Qco(X))$. As particular instances, we recover
recent results on the flat model
structure for quasi-coherent sheaves. Our approach also includes the case of
(infinite-dimensional) vector bundles, and of
restricted flat Mittag-Leffler quasi-coherent sheaves, as introduced by Drinfeld.
Finally, we prove that the unrestricted case does not induce a model category
structure as above in general.
\end{abstract}

\maketitle

\section{Introduction}

Let $X$ be a scheme and $\Qco(X)$ the category of all quasi-coherent sheaves on~$X$.
A convenient way of approaching the derived category $\mathbb D (\Qco (X))$ goes
back to
Quillen \cite{Quillen}, and consists in introducing a model category structure on
$\Ch (\Qco (X))$,
the category of unbounded chain complexes on $\Qco(X)$. In particular, one can compute
morphisms between two objects $X$ and $Y$ of $\mathbb D (\Qco (X))$ as the $\Ch
(\Qco (X))$-morphisms
between cofibrant and fibrant replacements of $X$ and $Y$, respectively, modulo
chain homotopy.

Recently, Hovey has shown that model category structures naturally arise from  small
cotorsion pairs over $\Ch (\Qco(X))$,
\cite{hovey2}. Since $\Qco(X)$ is a Grothendieck category \cite{EE}, there is a
canonical injective model category
structure on $\Ch (\Qco (X))$. However, this structure is not monoidal, that is,
compatible with the tensor product on $\Qco(X)$,
\cite[pp.\ 111-2]{Hov}. Another natural, but not monoidal, model structure on
$\Ch(\Qco(X))$ was constructed in \cite{Hov2}
under the assumption of $X$ being a~Noetherian separated scheme with enough locally
frees.

The lack of compatibility with the tensor product was partially solved in
\cite{G,Mur} by using flat quasi-coherent sheaves.
The main result of \cite{G} shows that in case $X$ is quasi-compact and
semi-separated, it is possible to construct a monoidal flat model structure
on $\Ch(\Qco(X))$. The weak equivalences of this model structure are the
same as the ones for the injective model structure, hence they induce
the same cohomology functors (see \cite{Mur} for a different approach).
However, the structure of flat quasi-coherent sheaves is rather
complex, and it is difficult to compute the associated fibrant and cofibrant
replacements. Moreover, the methods of \cite{G} depend heavily on the fact that
the class of all flat modules is closed under direct limits.

A different approach has recently been suggested in \cite{EEGR} for the particular
case of quasi-coherent sheaves on the
projective line ${\bf P^1}(k)$. In that paper it was shown that the class of
infinite-dimensional vector bundles
(i.e., those quasi-coherent sheaves whose sections in all open affine sets are
projective) imposes a monoidal model
category structure on $\Ch(\Qco({\bf P^1}(k))$. The proofs and techniques in
\cite{EEGR} are strongly based
on the Grothendieck decomposition theorem for vector bundles over the projective
line \cite{GRO}, hence they cannot be
extended to more general situations.

In the present paper, we show that the main results of \cite{EEGR} and
\cite{G} are particular instances of the following general theorem that
provides for a variety of model category structures on $\mathbb{C}(\Qco(X))$
parametrized by sets $S_v$ ($v \in V$) of modules of sections
(see Section \ref{qmodel} for unexplained terminology):

\begin{theorem}\label{maint}
Let $X$ be a semi-separated scheme. There is a model category structure on
$\mathbb{C}(\Qco(X))$ in which the
weak equivalences are the homology isomorphisms, the cofibrations (resp.\ trivial
cofibrations) are the monomorphism
with cokernels in $\dgC$ (resp.\ $\barC$), and the fibrations (resp.\ trivial
fibrations) are the epimorphisms whose kernels are
in $\dgCo$ (resp.\ $\barCo$). Moreover, if every $M\in S_v$ is a flat $\mathscr
R(v)$-module, and
$M\otimes_{\mathscr R(v)} N\in S_v$ for all $M, N\in S_v$, then the model category
structure is monoidal.
\end{theorem}
The proof of Theorem \ref{maint} is based on new tools for handling filtrations of
quasi--coherent sheaves developed in this paper. Thus it avoids the usual assumption
of closure under direct limits.

Theorem \ref{maint} immediately yields the following generalization of \cite[Theorem
6.1]{EEGR}:

\begin{corollary}\label{cor1}
Let $X$ be a scheme having enough infinite-dimensional vector
bundles (for example, a quasi-compact and quasi-separated scheme
that admits an ample family of invertible sheaves, or a noetherian,
integral, separated, and locally factorial scheme). Let
$\mathcal{C}$ be the class of all vector bundles on $X$.

Then there is a monoidal model category structure on $\mathbb{C}(\Qco(X))$ where
weak equivalences are
homology isomorphisms, the cofibrations (trivial cofibrations) are the monomorphisms
whose cokernels are $dg$-complexes
of vector bundles (exact complexes of vector bundles whose every cycle is a
vector bundle), and the fibrations
(trivial fibrations) are the epimorphisms whose kernels are in $\dgCo$ ($\barCo$).
\end{corollary}

Similarly, we immediately recover
\cite[Theorem 6.7]{G}:

\begin{corollary}\label{cor2}
Let $X$ be a scheme with enough flat quasi-coherent sheaves (for instance, let $X$
be quasi-compact and semi-separated,
see \cite[Proposition 16]{Mur}). Then there there is a monoidal model category
structure on $\Ch(\Qco(X))$ where weak
equivalences are homology isomorphisms, the cofibrations (trivial cofibrations) are
the monomorphisms whose cokernels
are dg-flat complexes (flat complexes). The fibrations (trivial fibrations) are the
epimorphisms whose kernels are
dg-cotorsion complexes (cotorsion complexes).
\end{corollary}

However, there are further interesting applications of Theorem \ref{maint}. Drinfeld
has
proposed quasi-coherent sheaves whose sections at affine open sets are flat and
Mittag-Leffler modules (in the sense of Raynaud and Gruson \cite{RG}) as the
appropriate objects
defining infinite-dimensional vector bundles on a scheme, see \cite[p.266]{D}. Here
we call
such quasi-coherent sheaves the Drinfeld vector bundles, and show that the
restricted ones
(bounded by a cardinal $\kappa$) fit into another instance of Theorem \ref{maint}:

\begin{corollary}\label{cor3}
Let $X$ be a semi-separated scheme possessing a generating set
$\mathcal G$ of Drinfeld vector bundles. Let $\kappa$ be an infinite
cardinal such that $\kappa  \geq |E|$ (in the notation of Section 3)
and each $\mathscr M \in \mathcal G$ is $\leq \kappa$-presented. For
each $v \in V$, let $S_v$ denote the class of all $\leq
\kappa$-presented flat Mittag-Leffler modules. Denote by $\mathcal
C$ the class of all Drinfeld vector bundles $\mathscr M$ such that
${\mathscr M}(v)$ has a $S_v$-filtration for each $v \in V$.

Then there is a monoidal model category structure on $\mathbb{C}(\Qco(X))$ where
weak equivalences are homology isomorphisms, the cofibrations (trivial cofibrations)
are monomorphism with cokernels in $\dgC$ ($\barC$), and the fibrations (trivial
fibrations)
are epimorphisms whose kernels are in $\dgCo$ ($\barCo$).
\end{corollary}

The reader may wonder whether it is possible to apply Theorem \ref{maint} to the
entire class
of Drinfeld vector bundles and impose thus a (monoidal) model category structure on
$\Ch(\Qco(X))$.
Our final theorem shows that this is not the case in general.
We adapt a recent consistency result of Eklof and Shelah \cite{ES} concerning
Whitehead groups to this setting, and
prove (in ZFC):

\begin{theorem}\label{notmodel}
The class $\mathcal D$ of all flat Mittag-Leffler abelian groups is not precovering.
Thus $\mathcal D$ cannot induce a cofibrantly generated model category structure on
$\Qco(\rm{Spec}(\mathbb{Z}))\cong \rm{Mod-}\mathbb{Z}$ compatible with its abelian
structure.
\end{theorem}

\section{Notation and Preliminaries}

Let ${\mathcal A}$ be a
Grothendieck category. A well-ordered direct system of objects of~$\mathcal A$, $ (
A_{ \alpha } \mid \alpha
\leq \lambda ) $, is said to be \emph{continuous} if $ A_0  = 0 $ and, for
each limit ordinal $ \beta \leq \lambda $, we have $ \displaystyle
A_{\beta}=\lim_{ \rightarrow } \:  A_{ \alpha } \:$ where the limit is taken over
all ordinals $\alpha{<}\beta$. A~continuous direct system $ ( A_{ \alpha}
\mid\alpha \leq \lambda ) $ is called a \emph{continuous directed union} if
all morphisms in the system are monomorphisms.

\begin{definition}\label{defilt}
Let $ {\mathcal L} $ be a class of objects of ${\mathcal A}$. An
object $A$ of $\mathcal A$ is \emph{$\mathcal L$-filtered} if  $ A =
\displaystyle \lim_{ \rightarrow } \: A_{ \alpha } $ for a
continuous directed union $ ( A_{ \alpha } \mid  \alpha \leq \lambda
) $ satisfying that,  for each $ \alpha + 1 \leq \lambda
$, $ \mbox{Coker} \: ( A_{ \alpha } \rightarrow A_{ \alpha + 1 } )$
is isomorphic to an element of $ {\mathcal L} $.

We denote by $Filt(\mathcal L)$ the class of
all $\mathcal L$-filtered objects in $\mathcal A$. A class $ {\mathcal L} $ is
said to be \emph{closed under $\mathcal L$-filtrations} in case $Filt(\mathcal L)=
\mathcal L$.
\end{definition}

\begin{definition}
Let $\mathcal D$ be a class of objects of ${\mathcal{A}}$. We will denote by
$\mathcal D^{\perp}$ the subclass
of $\mathcal{A}$ defined by $$\mathcal D^\perp=\hbox{Ker}\Ext 1{\mathcal A}
{\mathcal D}{-} =\{Y\in {\mathcal Ob}({\mathcal A})\mid \Ext 1 {\mathcal A} D Y=0,\
{\rm{for\ all}}\
D\in {\mathcal D}\}.$$ Similarly,
$$^{\perp}{\mathcal D}=\hbox{Ker}\Ext 1{\mathcal A}{-}
{\mathcal D} =\{Z\in {\mathcal Ob}({\mathcal A})\mid\Ext 1 {\mathcal A} Z D=0,\
{\rm{for\ all}}\
D\in {\mathcal D}\}.$$
Analogously, we will define
$${\mathcal D}^{{\perp}_{\infty}}=\{Y\in {\mathcal Ob}({\mathcal A})\mid\Ext i
{\mathcal A} D Y=0,\ {\rm{for\ all}}\ D\in {\mathcal D}\ {\rm and}\ i\geq 1\} $$
and
$$^{\perp_{\infty}}{\mathcal D}=\{Z\in {\mathcal Ob}({\mathcal A})\mid\Ext i
{\mathcal A} Z D=0,\ {\rm{for\ all}}\
D\in {\mathcal D}\ {\rm and}\ i\geq 1\}. $$
\end{definition}
Let us recall the following definitions from \cite{GT}.
\begin{definition}
A pair (${\mathcal F}, {\mathcal C}$) of classes of objects of
${{\mathcal{A}}}$ is called a \emph{cotorsion pair} if ${\mathcal
F}^{\perp}={\mathcal C}$ and if $^{\perp}{ \mathcal C}={\mathcal F}$. The cotorsion
pair is said
to have \emph{enough injectives} (resp.\ \emph{enough projectives}) if for each
object $Y$
of ${{\mathcal{A}}}$ there exists an exact sequence
 $0\rightarrow Y\rightarrow C\rightarrow
F\rightarrow 0$ (resp.\ for each object $Z$ of $\mathcal A$ there exists an
exact sequence $0\rightarrow C'\rightarrow F'\rightarrow
Z\rightarrow 0$) such that $F,F'\in {\mathcal F}$ and $C,C'\in {\mathcal
C}$. A cotorsion pair $({\mathcal F},{\mathcal C})$ is \emph{complete}
provided it has enough injectives and enough projectives.

\end{definition}

The proof of the following lemma is the same as for module categories (see
\cite[Lemma 2.2.10]{GT}).
\begin{lemma}
Let $\mathcal A$ be a Grothendieck category with enough projectives and
let $(\mathcal F,\mathcal C)$ be a cotorsion pair on $\mathcal A$. The following
conditions are equivalent
\medskip\par\noindent
$\rm a)$ If $0\to F'\to F\to F''\to 0$ is exact with $F,F''\in
\mathcal F$, then $F'\in\mathcal F$.
\medskip\par\noindent
$\rm b)$ If $0\to C'\to C\to C''\to 0$ is exact with $C',C\in \mathcal C$,
then $C''\in\mathcal C$.
\medskip\par\noindent
$\rm c)$ ${\rm Ext}^2(F,C)=0$ for all $F\in \mathcal F$ and $C\in \mathcal C$.
\medskip\par\noindent
$\rm d)$ ${\rm Ext}^n(F,C)=0$ for all $n\geq 1$ and all $F\in \mathcal F$ and
$C\in \mathcal C$.
\end{lemma}

A cotorsion pair satisfying the equivalent conditions above is
called \emph{hereditary}. So $(\mathcal F,\mathcal C)$ is a hereditary cotorsion
pair, if and only if $\mathcal F={}^{\perp_{\infty}}{\mathcal C}$ and $\mathcal
C={\mathcal F}^{{\perp}_{\infty}}$.

We finish this section by recalling the notion of a Mittag-Leffler module from
\cite{RG}.

\begin{definition}\label{flml} \rm Let $R$ be a ring and $M$ a right $R$-module.
Then $M$ is {\em Mittag--Leffler} provided that the canonical map
$M \otimes_R \prod_{i \in I} M_i \to \prod_{i \in I} M \otimes_R M_i$ is
monic for each family of left $R$--modules $( M_i \mid i \in I )$.
\end{definition}

For example, all finitely presented modules, and all projective modules, are
Mittag--Leffler.
Any countably generated flat Mittag-Leffler module is projective. In fact,
projectivity of a module $M$ is
equivalent to $M$ being flat Mittag-Leffler and a direct sum of countably generated
submodules (see \cite{RG} and \cite[Theorem 2.2]{D}).

We refer to \cite{EM,GT,HAR,Hov} for unexplained terminology used in this paper.

\section{Filtrations of Quasi--Coherent Sheaves}\label{filtqc}

Let $X$ be a scheme. Let $Q_X=(V,E)$ be the quiver whose set, $V$,
of vertices is a subfamily of the family of all open affine sets of
$X$ such that $V$ covers both $X$ and all intersections $O\cap O'$
of open affine sets $O,O'$ of $V$. The set of edges, $E$, consists
of the reversed arrows $v\to u$ corresponding to the inclusions
$u\subseteq v$ where $u$ and $v$ are in $V$. We say that $Q_X$ is a
{\em quiver associated to the scheme} $X$. Note that this quiver is
not unique, because different choices of the set of vertices $V$ may
give rise to non-isomorphic quivers associated to the same scheme
$X$.

As explained in \cite[Section 2]{EE},
there is an equivalence between the category of quasi--coherent sheaves on $X$
and the category of quasi--coherent $\mathscr R$--modules where
$\mathscr R$ is the representation of the quiver $Q_X$ by the sections of the
structure
sheaf $\mathcal{O}_X$. A~quasi--coherent sheaf $\mathcal F$ on $X$ corresponds to a
\emph{quasi--coherent $\mathscr R$-module} $\mathscr M$ defined by the following data:

\begin{enumerate}
\item An {\em $\mathscr R$-module} on $X$, that is, an ${\mathscr R}(u)$-module
$M_u$, for each
$u\in V$ and a ${\mathscr R}(u)$--morphism $\rho_{uv}:M_u\to M_v$ for each edge $u\to
v$ in $E$;
\item The {\em quasi--coherence} condition, saying that the induced morphism
$$id_{{\mathscr R}(v)}\otimes \rho_{uv}: {\mathscr R}(v)\otimes_{{\mathscr
R}(u)}M_u\to
{\mathscr R}(v)\otimes_{{\mathscr R}(u)}M_v\cong M_v$$  is an $\mathscr
R(v)$-isomorphism, for each arrow $u\to v$ in $E$;
\item The {\em compatibility condition}, saying that if $w\subseteq v\subseteq u$,
with $w,v,u\in V$, then $\rho_{uw} = \rho_{vw} \circ \rho_{uv}$.
\end{enumerate}
Note that quasi--coherent subsheaves $\mathscr F'$ of $\mathscr F$ correspond to
quasi--coherent $\mathscr R$-sub\-mo\-du\-les $\mathscr M'$ of $\mathscr M$ (where
the latter means that $M'_v$ is an $\mathscr R(v)$--submodule of $M_v$ and the map
$\rho'_{uv}$ is a restriction of $\rho_{uv}$ for each edge $u\to v$ in $E$). If
$(\mathscr F_i)_{i\in I}$ are quasi--coherent subsheaves of $\mathscr F$ then
$\mathscr F'=\sum_{i\in I}\mathscr F_i$ (resp. $\mathscr F'=\mathscr F_1\cap
\mathscr F_2$) corresponds to the quasi-coherent submodule $\mathscr M'$ such that
$M'_v=\sum_{i\in I}(M_i)_v$ (resp. such that $M'_v=(M_1)_v\cap (M_2)_v$) and the
maps $\rho'_{uv}$ are restrictions of $\rho_{uv}$.

Recall that $\Qco(X)$ denotes the category of all quasi--coherent sheaves on $X$.
This is a Grothendieck category by
\cite[p.290]{EE}. Note that in our setting, if $u \subseteq v$ are affine open
subsets in $V$,
then ${\mathscr R}(u)$ is a flat ${\mathscr R}(v)$--module, see \cite[III.9]{HAR}.

Recall that a quasi-coherent sheaf $\mathscr M$ on $X$ is a (classical algebraic)
vector bundle if ${\mathscr M}(u)$ is a free ${\mathscr R}(u)$-module of finite rank
for every open affine set $u$. In this paper we adopt the following more general
definition: $\mathscr M$ is a {\em vector bundle} if ${\mathscr M}(u)$
is a (not necessarily finitely generated) projective ${\mathscr R}(u)$--module for
each open affine set $u$ (see \cite[\S2. Definition]{D}).

In \cite[Section 2.Remarks]{D}, Drinfeld proposed to consider the following more
general notion of a vector bundle
(see also \cite[Appendices 5 and 6]{D2}). Thus, we call a quasi--coherent sheaf
$\mathscr M$ a \emph{Drinfeld vector bundle}
provided that ${\mathscr M}(u)$ is a flat Mittag--Leffler ${\mathscr R}(u)$--module
for each open affine set $u$ (cf.\ \cite[p.266]{D}).

\medskip
One of the main goals of this paper is to construct monoidal model category
structures associated to these generalized notions of vector bundles. In order to
achieve this aim we will need to characterize these classes as closures under
filtrations of certain of their subsets.

The following tools will play a central role in our study of these filtrations, both
in the case of modules over a ring, and of quasi-coherent sheaves on a scheme.

The first tool is known as Eklof's Lemma (see \cite[Theorem 1.2]{E}):

\begin{lemma}\label{eklof}
Let $R$ be a ring and $\mathcal C$ be a class of modules. Let $M$ be a module
possesing a ${}^\perp \mathcal C$-filtration.
Then $M \in {}^\perp \mathcal C$.
\end{lemma}

\begin{remark}\label{remek}
The proof of Lemma \ref{eklof} given in \cite[Lemma 3.1.2]{GT} needs only
embeddability of each module into an injective one, so the lemma holds in $\Qco
(X)$, and in fact in any Grothendieck category.
\end{remark}

Our second tool is known as Hill's Lemma
(see \cite[Theorem 4.2.6]{GT}, \cite[Lemma 1.4]{SaTr}, or \cite[Theorem 6]{ST}).
It will allow us to extend a given filtration of a module $M$ to a complete lattice
of its submodules having similar properties.

\begin{lemma}\label{hill}
Let $R$ be a ring, $\lambda$ a regular infinite cardinal, and $\mathcal J$ a class
of ${<} \lambda$--presented modules.
Let $M$ be a module with a $\mathcal J$--filtration $\mathcal M = (M_\alpha \mid
\alpha \leq \sigma )$.
Then there is a family $\mathcal H$ consisting of submodules of $M$ such that
\begin{enumerate}
\item $\mathcal M \subseteq \mathcal H$,
\item $\mathcal H$ is closed under arbitrary sums and intersections,
\item $P/N$ has a $\mathcal J$--filtration for all $N, P \in \mathcal H$ such that $N
\subseteq P$, and
\item If $N \in \mathcal H$ and $T$ is a subset of $M$ of cardinality ${<} \lambda$,
then there exists $P \in \mathcal H$ such that $N \cup T \subseteq P$ and $P/N$ is
${<} \lambda$--presented.
\end{enumerate}
\end{lemma}

We will also need the following application of Lemma \ref{hill} (see \cite[Theorem
4.2.11]{GT} and \cite[Theorem 10]{ST}):

\begin{lemma}\label{cotorp} Let $R$ be a ring, $\lambda$ a regular uncountable
cardinal,
and $\mathcal J$ a class of ${<} \lambda$--presented modules. Let $\mathcal A =
{}^\perp (\mathcal J ^\perp)$, and let
$\mathcal A ^{{<}\lambda} $ denote the class of all ${<} \lambda$--presented modules from
$\mathcal A$.
Then every module in $\mathcal A$ is $\mathcal A^{{<}\lambda}$--filtered.
\end{lemma}

If $\kappa$ is a cardinal and $\mathscr M$ a quasi--coherent sheaf, then $\mathscr M$
is called {\em locally $\leq \kappa$--presented} if for each $v \in V$, the
${\mathscr R}(v)$--module $\mathscr M (v)$ is $\leq \kappa$--presented.
Notice that if $\kappa \geq|V|$ and $\kappa\geq|{\mathscr R}(v)|$ for each $v \in V$,
then this definition is equivalent to saying that $\mathscr M$ is
$\kappa^+$--presentable
in the sense of \cite[Lemma 6.1]{G}, and also to $|\bigoplus_{v \in V} M(v)| \leq
\kappa$.

\medskip
For future use in Section \ref{qmodel} we now present a version of Hill's Lemma for
the category $\Qco(X)$.
For this version, we assume that  $X$ is a scheme, $\lambda$ a regular infinite
cardinal such that $\lambda {>} |V|$ and
$\lambda {>} |{\mathscr R}(v)|$ for all $v \in V$, and $\mathcal J$ a class of locally
${<}\lambda$--presented objects of $\Qco(X)$.
Further, let $\mathscr M$ be a quasi--coherent sheaf possessing  a $\mathcal
J$--filtration $\mathcal O = (\mathscr M_\alpha \mid \alpha \leq \sigma )$.

By \cite[Corollary 2.3]{EEGR0}, there exist locally ${<}\lambda$--presented
quasi-coherent sheaves $\mathscr A_\alpha\subseteq \mathscr M_{\alpha+1}$ such that
$\mathscr M_{\alpha +1}=\mathscr M_\alpha + \mathscr A_\alpha$ for each $\alpha {<}
\sigma$. A set $S\subseteq \sigma$ is called \emph{closed} provided that $\mathscr
M_\alpha\cap \mathscr A_\alpha \subseteq \sum_{\beta{<}\alpha, \beta\in S}\mathscr
A_{\beta}$ for each $\alpha\in S$.

\begin{lemma}\label{HillQco} Let $\mathcal H=\{\sum_{\alpha\in S}\mathscr A_\alpha
\mid S\textrm{ closed}\,\}$.
Then $\mathcal H$ satisfies the following conditions:
\begin{enumerate}
\item $\mathcal O \subseteq \mathcal H$,
\item $\mathcal H$ is closed under arbitrary sums,
\item $\mathscr P/\mathscr N$ has a $\mathcal J$--filtration whenever $\mathscr
N,\mathscr P \in \mathcal H$ are such that $\mathscr N
\subseteq \mathscr P$.
\item If $\mathscr N \in \mathcal H$ and $\mathscr X$ is a locally ${<}
\lambda$--presented quasi--coherent subsheaf of $\mathscr M$,
then there exists $\mathscr P \in \mathcal H$ such that $\mathscr N + \mathscr X
\subseteq \mathscr P$ and $\mathscr P/\mathscr N$ is locally
${<} \lambda$--presented.
\end{enumerate}
\end{lemma}

\noindent\Proof Note that for each ordinal $\alpha\leq \sigma$, we have $\mathscr
M_\alpha=\sum_{\beta{<}\alpha} \mathscr A_{\beta}$, hence $\alpha$ is a closed subset
of $\sigma$. This proves condition $(1)$. Since any union of closed subsets is
closed, condition $(2)$ holds.

In order to prove condition $(3)$, we consider closed subsets $S,T$ of $\sigma$ such
that $\mathscr N=\sum_{\alpha\in S}\mathscr A_{\alpha}$ and $\mathscr
P=\sum_{\alpha\in T}\mathscr A_{\alpha}$. Since $S\cup T$ is closed, we will
w.l.o.g.\ assume that $S\subseteq T$. We define a $\mathcal J$--filtration of
$\mathscr P/\mathscr N$ as follows. For each $\beta \leq \sigma$, let $\mathscr
F_{\beta}=(\sum_{\alpha\in T\setminus S,\alpha{<}\beta}\mathscr A_{\alpha}+\mathscr
N)/\mathscr N$. Then $\mathscr F_{\beta+1}=\mathscr F_{\beta}+(\mathscr
A_{\beta}+\mathscr N)/\mathscr N$ for $\beta\in T\setminus S$ and $\mathscr
F_{\beta+1}=\mathscr F_{\beta}$ otherwise.

Let $\beta\in T\setminus S$. Then $ \mathscr F_{\beta+1}/\mathscr F_{\beta}\cong
\mathscr A_{\beta}/(\mathscr A_{\beta}\cap(\sum_{\alpha\in T\setminus
S,\alpha{<}\beta}\mathscr A_{\alpha}+\mathscr N)),$ and since $\beta\in T\setminus S$
and $T$ is closed, we have
$$\mathscr A_{\beta}\cap(\sum_{\alpha\in T\setminus S,\alpha{<}\beta}\mathscr
A_{\alpha}+\mathscr N)=\mathscr A_{\beta}\cap (\sum_{\alpha\in
S,\alpha{>}\beta}\mathscr A_{\alpha}+\sum_{\alpha\in T,\alpha{<}\beta}\mathscr
A_{\alpha})\supseteq$$ $$\supseteq \mathscr A_{\beta}\cap (\sum_{\alpha\in S,\alpha{>}
\beta}\mathscr A_{\alpha}+(\mathscr M_{\beta}\cap \mathscr A_{\beta}))\supseteq
\mathscr M_{\beta}\cap \mathscr A_{\beta}.$$
Let $\mathscr B_{\beta}=\sum_{\alpha\in S,\alpha{>}\beta}\mathscr
A_{\alpha}+\sum_{\alpha\in T,\alpha{<}\beta}\mathscr A_{\alpha}$
We will prove that $\mathscr A_{\beta}\cap\mathscr B_{\beta}=\mathscr M_{\beta}\cap
\mathscr A_{\beta} $. We have only to show that for each $v\in V$, $\mathscr
A_{\beta}(v)\cap \mathscr B_{\beta}(v)\subseteq \mathscr A_{\beta}(v)\cap \mathscr
M_{\beta}(v)$. Let $a\in \mathscr A_{\beta}(v)\cap \mathscr B_{\beta}(v)$. Then
$a=c+a_{\alpha_0}+\cdots +a_{\alpha_k}$ where $c\in \sum_{\alpha\in
T,\alpha{<}\beta}\mathscr A_{\alpha}(v)\subseteq \mathscr M_{\beta}(v)$, $\alpha_i\in
S$ and $a_{\alpha_i}\in \mathscr A_{\alpha_i}(v)$ for all $i\leq k$ and
$\alpha_i{>}\alpha_{i+1}$ for all $i{<}k$. W.l.o.g., we can assume that $\alpha_0$ is
minimal possible. If $\alpha_0{>}\beta$, then $a_{\alpha_0}=a-c-a_{\alpha_1}+\cdots
-a_{\alpha_k}\in \mathscr M_{\alpha_0}(v)\cap \mathscr A_{\alpha_0}(v)\subseteq
\sum_{\alpha\in S,\alpha{<}\alpha_0} \mathscr A_{\alpha}(v)$ (since $\alpha_0\in S$),
in contradiction with the minimality of $\alpha_0$. Since $\beta \notin S$, we infer
that $\alpha_0{<}\beta$, $a\in \mathscr M_{\beta}(v)$, and $\mathscr A_{\beta}\cap
\mathscr B_{\beta}=\mathscr A_{\beta}\cap \mathscr M_{\beta}$.

So if $\beta \in T\setminus S$ then $\mathscr F_{\beta+1}/\mathscr F_{\beta}\cong
\mathscr A_{\beta}/(\mathscr M_{\beta}\cap \mathscr A_{\beta})\cong \mathscr
M_{\beta+1}/\mathscr M_{\beta}$, and the latter is isomorphic to an element of
$\mathcal J$ because $\mathcal O$ is a $\mathcal J$--filtration of $\mathscr M$.
This finishes the proof of condition $(3)$.

For condition $(4)$ we first claim that each subset of $\sigma$ of
cardinality ${<}\lambda$ is contained in a closed subset of
cardinality ${<}\lambda$. Since $\lambda$ is regular and unions of
closed sets are closed, it suffices to prove the claim only for
one--element subsets of $\sigma$. By induction on $\beta$ we prove
that each $\beta{<}\sigma$ is contained in a closed set $S$ of
cardinality ${<}\lambda$. If $\beta{<}\lambda$ we take $S=\beta+1$.

Otherwise, consider the short exact sequence $0\to \mathscr
M_{\beta}\cap \mathscr A_{\beta}\to \mathscr A_{\beta}\to \mathscr
M_{\beta+1}/\mathscr M_{\beta}\to 0 $. By our assumption on
$\lambda$, since $A_{\beta}$ is locally ${<}\lambda$--presented, so is
$\mathscr M_{\beta}\cap \mathscr A_{\beta}$. Hence for each $v \in
V$, $\mathscr M_{\beta}(v) \cap \mathscr A_{\beta}(v) \subseteq
\sum_{\alpha \in S_v}\mathscr A_{\alpha}(v)$ for a subset $S_v
\subseteq \beta$ of cardinality ${<} \lambda$. By our inductive
premise, the set $\bigcup_{v \in V} S_v$ is contained in a closed
subset $S^\prime$ of cardinality ${<} \lambda$. Let $S=S^\prime \cup
\{\beta\}$. Then $S$ is closed because $S^\prime$ is closed, and
$\mathscr M_{\beta}\cap \mathscr A_{\beta}\subseteq \sum_{\alpha \in
S^\prime}\mathscr A_{\alpha}$.

Finally if $\mathscr N=\sum_{\alpha\in S}\mathscr A_{\alpha}$ and
$\mathscr X$ is a locally ${<}\lambda$--presented quasi--coherent
subsheaf of $\mathscr M$, then $\mathscr X\subseteq \sum_{\alpha\in
T}\mathscr A_{\alpha}$ for a subset $T$ of $\sigma$ of cardinality
${<}\lambda$. By the above we can assume that $T$ is closed and put
$\mathscr P=\sum_{\alpha \in S \cup T}\mathscr A_{\alpha}$. By (the
proof of) condition $(3)$ $\mathscr P/\mathscr N$ is $\mathcal
J$--filtered, and the length of the filtration can be taken $\leq
|T\setminus S|{<}\lambda$. This implies that $\mathscr P/\mathscr N$
is locally ${<}\lambda$--presented. \qed

\medskip
Our third tool is essentially \cite[Proposition 3.3]{EE} (where we omit the
condition of ${\mathscr M}'(v)$ being a pure submodule in ${\mathscr M}(v)$,
because we do not need it in the sequel). This tool will be applied to form filtrations
of quasi--coherent sheaves by connecting the individual ${\mathscr R}(v)$--module
filtrations
for all $v\in V$.

\begin{lemma}\label{nocont}
Let $Q_X = (V,E)$ be a quiver associated to a scheme $X$, and let $\mathscr M\in
\Qco(X)$.
Let $\kappa$ be an infinite cardinal such that $\kappa\geq|V|$, and
$\kappa\geq|{\mathscr R}(v)|$ for all $v\in V$. Let $X_v\subseteq {\mathscr M}(v)$
be subsets with $|X_v|\leq \kappa$ for all
$v \in V$. Then there is a locally $\leq \kappa$-presented quasi-coherent subsheaf
${\mathscr M}'\subseteq {\mathscr M}$ such that $X_v\subseteq {\mathscr M}'(v)$ for
all $v \in V$.
\end{lemma}

Now we fix our notation:

\begin{notation}\label{nttn} Let $Q_X =(V,E)$ be a quiver associated to a scheme $X$,
and $\kappa$ be an infinite cardinal such that $\kappa\geq|V|$ and
$\kappa\geq|{\mathscr R}(v)|$ for all $v\in V$. For each $v \in V$,
let $S_v$ be a class of $\leq \kappa$--presented ${\mathscr
R}(v)$--modules, $\mathcal F _v = {}^\perp(S_v ^\perp)$, $\mathcal
L$ be the class of all locally $\leq \kappa$--presented
quasi--coherent sheaves $\mathscr N$ such that $\mathscr N(v) \in
\mathcal F _v$ for each $v\in V$, and $\mathcal C$ be the class of
all quasi--coherent sheaves $\mathscr M$ such that ${\mathscr
M}(v)\in {\mathcal F}_v$ for each $v \in V$.
\end{notation}

\begin{theorem}\label{prevnew}
Each quasi--coherent sheaf $\mathscr M\in \mathcal C$ has an $\mathcal L$--filtration.
\end{theorem}
\noindent\Proof  Let $v \in V$ and put $\lambda = \kappa ^+$. Denote by $\mathcal F
_v ^{\leq\kappa}$ the subclass of $\mathcal F _v$ consisting of all $\leq
\kappa$-presented modules. By Lemma \ref{cotorp}, ${\mathscr M}(v)$ has a $\mathcal
F _v ^{\leq\kappa}$-filtration $\mathcal M _v$. Denote by $\mathcal H _v$ the family
associated to $\mathcal M _v$ in Lemma \ref{hill}. And let $\{ m_{v,\alpha} \mid
\alpha {<} \tau_v \}$ be an ${\mathscr R}(v)$-generating set of the ${\mathscr
R}(v)$-module ${\mathscr M}(v)$. W.l.o.g., we can assume that $\tau = \tau _v$ for
all $v \in V$.

We will construct an $\mathcal L$--filtration $({\mathscr M}_\alpha \mid \alpha \leq
\tau )$ of $\mathscr M$ by induction on $\alpha$.
Let ${\mathscr M}_0 = 0$. Assume that ${\mathscr M}_\alpha$ is defined for some
$\alpha {<} \tau$ so that
$ {\mathscr M}_{\alpha}(v) \in \mathcal H _v$ and $m_{v,\beta} \in {\mathscr M}(v)$
for all $\beta {<} \alpha$ and all $v \in V$.
Set $N_{v,0} = {\mathscr M}_{\alpha}(v)$. By Lemma \ref{hill}.(4), there is a module
$N_{v,1} \in \mathcal H _v$ such that $N_{v,0} \subseteq N_{v,1}$, $m_{v,\alpha} \in
N_{v,1}$ and $N_{v,1}/N_{v,0}$ is $\leq \kappa$-presented.

By Lemma \ref{nocont} (with ${\mathscr M}$ replaced by ${\mathscr M}/{\mathscr
M}_\alpha$, and $X_v = N_{v,1}/{\mathscr M}_{\alpha}(v)$) there is a quasi-coherent
subsheaf ${\mathscr M}^{\prime}_1$ of ${\mathscr M}$ such that ${\mathscr M}_\alpha
\subseteq {\mathscr M}'_1$ and ${\mathscr M}'_1/{\mathscr M}_\alpha$ is locally
$\leq \kappa$-presented. Then ${\mathscr M}'_1(v)=N_{v,1} + \langle T_v \rangle$ for
a subset $T_v \subseteq {\mathscr M}'_1(v)$ of cardinality $\leq \kappa$, for each
$v \in V$.

By Lemma \ref{hill}.(4) there is a module $N_{v,2} \in \mathcal H _v$ such that
${\mathscr M}'_1(v)= N_{v,1} + \langle T \rangle \subseteq N_{v,2}$ and
$N_{v,2}/N_{v,1}$ is $\leq \kappa$-presented.

Proceeding similarly, we obtain a countable chain $({\mathscr M}'_n \mid n{<}\aleph_0)
$ of quasi-coherent subsheaves of ${\mathscr M}$, as well as a countable chain
$(N_{v,n} \mid n {<} \aleph_0 )$ of ${\mathscr R}(v)$-submodules of ${\mathscr M}(v)$,
for each $v \in V$. Let ${\mathscr M}_{\alpha + 1} =   \bigcup_{n {<} \aleph_0}
{\mathscr M}'_n$. Then ${\mathscr M}_{\alpha+1}$ is a quasi-coherent subsheaf of
${\mathscr M}$ satisfying ${\mathscr M}_{\alpha + 1}(v) = \bigcup_{n {<} \aleph_0}
N_{v,n}$ for each $v\in V$. By Lemma \ref{hill}.(2) and (3) we deduce that
${\mathscr M}_{\alpha + 1}(v) \in \mathcal H _v$ and ${\mathscr M}_{\alpha +
1}(v)/{\mathscr M}_\alpha (v) \in \mathcal F _v ^{\leq\kappa}$. Therefore ${\mathscr
M}_{\alpha + 1}/{\mathscr M}_\alpha \in \mathcal L$.

Assume ${\mathscr M}_\beta$ has been defined for all $\beta {<} \alpha$ where $\alpha$
is a limit ordinal $\leq \tau$. Then we define ${\mathscr M}_\alpha = \bigcup_{\beta
{<} \alpha} {\mathscr M}_\beta$. Since $m_{v,\alpha} \in {\mathscr M}_{\alpha +1}(v)$
for all $v \in V$ and $\alpha {<} \tau$, we have ${\mathscr M}_\tau(v) = {\mathscr
M}(v)$, so
$({\mathscr M}_\alpha \mid \alpha \leq \tau )$ is an $\mathcal L$--filtration of
${\mathscr M}$. \qed

\begin{remark}\label{remk}
Recall that a module $N$ is {\em strongly $\leq \kappa$-presented} provided that $N$
has a projective resolution consisting of $\leq \kappa$-generated projective
modules. If this is the case we will always consider only the projective resolutions
of $N$ that consist of $\leq \kappa$-generated modules.

A class of modules $\mathcal C$ is {\em syzygy closed} if for each $C \in \mathcal
C$, the first (and hence each) syzygy of $C$ in some projective resolution of $C$ is
contained in $\mathcal C$.

We note that Theorem \ref{prevnew} remains true under the stronger assumption that
for each $v \in V$, $S_v$ is a class of strongly $\leq \kappa$-presented
${\mathscr R}(v)$-modules and $\mathcal F _v = {}^{\perp_\infty}(S_v
^{\perp_\infty})$.
\end{remark}

It is clear that the class $\mathcal C$ is closed under extensions,
retractions and direct sums.

As a consequence of Theorem \ref{prevnew} we get the following two
corollaries.

\begin{corollary}\label{cfilt}
Let $X$ be any scheme with associated quiver $Q_X$. Let $\mathcal C$ and $\mathcal
L$ be the subclasses of $\Qco(X)$ defined above. Then $\mathcal C={Filt}(\mathcal L)$.
\end{corollary}
\noindent\Proof The inclusion $\mathcal C\subseteq {Filt}(\mathcal L)$ follows by
Theorem \ref{prevnew},
and ${Filt}(\mathcal L) \subseteq {Filt}(\mathcal C) \subseteq \mathcal C$ by Lemma
\ref{eklof}
(and Remark \ref{remek}). \qed

\begin{corollary}\label{und}
Let $X$ be any scheme with associated quiver $Q_X$. Let $\mathcal C$ and $\mathcal
L$ be the subclasses of $\Qco(X)$ defined above. Suppose that $\mathcal C$ contains
a generator of $\Qco(X)$. Then $(\mathcal C,{\mathcal L} ^{\perp})$ is a complete
cotorsion pair.
\end{corollary}
\noindent\Proof Since $\mathcal L\subseteq {}^{\perp}({\mathcal L} ^{\perp})$, we
have ${Filt}(\mathcal L)\subseteq {Filt}({}^{\perp}({\mathcal L}
^{\perp})) $. By Lemma \ref{eklof}, ${Filt}({}^{\perp}({\mathcal L}
^{\perp}))={}^{\perp}({\mathcal L} ^{\perp}) $. So by Corollary \ref{cfilt},
$\mathcal C\subseteq {}^{\perp}({\mathcal L} ^{\perp})  $.

In order to prove that $(\mathcal C,{\mathcal L} ^{\perp})$ is a complete cotorsion
pair, we first show that ${}^{\perp}({\mathcal L} ^{\perp}) \subseteq \mathcal C $.
By \cite[Lemma 2.4, Theorem 2.5]{EEGO}, for all $\mathscr Q\in \Qco(X)$ there exists a
short exact sequence\begin{equation}\label{equat1}
0\to \mathscr Q\to \mathscr P\to \mathscr Z\to 0
\end{equation}
where $\mathscr P\in \Uo$ and
$\mathscr Z$ has an $\mathcal L$-filtration. Given any $\mathscr M\in
\Qco(X)$, since the generator $\mathscr G$ of $\Qco(X)$ is in $\mathcal C$, there
exists a short
exact sequence
$$0\to \mathscr U\to \mathscr G'\to \mathscr M\to 0 $$ where $\mathscr G'$ is a
direct sum of copies of $\mathscr G\in {\mathcal L}$.
Now let
$$0\to \mathscr U\to {\mathscr N}\to \mathscr Z\to 0$$ be exact with ${\mathscr
N}\in \Uo$ and $\mathscr Z$ admitting an $\mathcal L$-filtration. Form a pushout and
get

$$\begin{CD}
@.     0@.        0      @.     @.\\
@.     @VVV       @VVV   @.     @.\\
0@>>>  {\mathscr U}@>>> {\mathscr G ^\prime}@>>>  {\mathscr M}@>>>  0\\
@.     @VVV  @VVV   @|     @.\\
0@>>>  {\mathscr N}@>>> {\mathscr W}@>>>  {\mathscr M}@>>>  0\\
@.     @VVV       @VVV @.  @.\\
@.     {\mathscr Z} @= {\mathscr Z} @.  @.\\
@.     @VVV       @VVV   @.     @.\\
@.     0@.        0      @.     @.
\end{CD}$$

Then since ${\mathscr G'}$ is a direct sum of copies of $\mathscr G\in {\mathcal C}$
and $\mathscr Z$ has an $\mathcal L$-filtration (so $\mathscr Z\in\mathcal C$ by
Corollary \ref{cfilt}), we see that $\mathscr W\in \mathcal C$. Also ${\mathscr
N}\in \Uo$.
Hence if ${\mathscr M}\in \oUo$ we get that $0\to {\mathscr N}\to \mathscr W\to
{\mathscr M}\to 0$ splits and
so ${\mathscr M}$ is a direct summand of $\mathscr W\in \mathcal C$. But then
$\mathscr M\in \mathcal C$ because $\mathcal C$ is closed under direct summands.

This proves that $\mathcal C= {}^{\perp}({\mathcal L} ^{\perp})$. Moreover
(\ref{equat1}) shows that the cotorsion pair $(\mathcal C,\mathcal L^\perp)$ has
enough injectives, and the second line of the diagram above that it has enough
projectives. \qed

\medskip
Focussing on particular classes of modules, we obtain several interesting
corollaries of Theorem \ref{prevnew}:

\begin{corollary}[Kaplansky Theorem for vector bundles]
Let $X$ be a scheme and $\kappa$ an infinite cardinal such that $\kappa\geq|V|$ and
$\kappa\geq|{\mathscr R}(v)|$ for all $v\in V$.
Then every vector bundle on $X$ has an $\mathcal L$--filtration
where $\mathcal L$ is the class of all locally $\leq \kappa$--presented
vector bundles.

In particular, if $X$ is a scheme, $Q_X = (V,E)$ is a quiver associated to $X$, and
both $V$ and all the rings ${\mathscr R}(v)$ ($v \in V$) are countable,
then every vector bundle on $X$ has a filtration by locally countably generated
vector bundles.
\end{corollary}
\noindent\Proof This follows by taking $S_v=\{{\mathscr R}(v)\}$ (so ${\mathcal
F}_v$ is the
class of all projective ${\mathscr R}(v)$-modules) for all $v\in V$, and then
applying Theorem \ref{prevnew}.
\qed

\medskip Let $\kappa$ be an infinite cardinal such that $\kappa \geq|V|$ and
$\kappa\geq|{\mathscr R}(v)|$ for all $v\in V$. For the next
corollary, we fix $\mathcal C$ to be the class of quasi--coherent
sheaves $\mathscr M$ such that each $\mathscr M(v)$ has a filtration
by $\leq \kappa$--presented flat Mittag--Leffler modules.

\begin{corollary}\label{boundrinf}
Let $X$ be a scheme, $\kappa$ an infinite cardinal such that $\kappa \geq|V|$ and
$\kappa\geq|{\mathscr R}(v)|$ for all $v\in V$. Let $\mathcal L$ be the
class of locally $\leq \kappa$-presented Drinfeld vector bundles.
For each $v \in V$, let $S_v$ denote the class of all $\leq \kappa$--presented flat
Mittag--Leffler modules. Let ${\mathcal F}_v$ and $\mathcal C$ be defined as
before Corollary \ref{cfilt}. Then ${Filt}(\mathcal L)=\mathcal C$.
\end{corollary}

In  \ref{nofilt}, we will see that in general Corollary \ref{boundrinf} fails for
arbitrary Drinfeld vector bundles. Our final application goes back to \cite[Section
4]{EE}:

\begin{corollary}
Let $X$ be a scheme. Let $\kappa$ be an infinite cardinal such that $\kappa \geq|V|$
and  $\kappa\geq|{\mathscr R}(v)|$ for all $v\in V$.

Then every flat quasi--coherent sheaf on $X$ has an $\mathcal L$--filtration
where $\mathcal L$ is the class of all locally $\leq \kappa$--presented flat
quasi--coherent sheaves.
\end{corollary}
\noindent\Proof For each vertex $v\in V$, we take a set $S_v$ of representatives of
iso classes of flat ${\mathscr R}(v)$--modules of cardinality $\leq \kappa$. Then by
Lemma \ref{eklof} and \cite[Lemma 1]{BBE} it follows that ${\mathcal
F}_v={}^\perp(S_v ^\perp)$ is the class of all flat ${\mathscr R}(v)$--modules.
Finally, we apply Theorem \ref{prevnew}.
\qed

\section{Quillen Model Category Structures on $\Ch(\Qco(X))$.}\label{qmodel}

In this section we develop a method for constructing a model structure on
${\Ch}(\Qco(X))$ starting from {\it a priori} given sets of modules over
sections of the structure sheaf associated to $X$. Our main tool will be Hovey's
Theorem relating cotorsion pairs to model category structures (see \cite[Theorem
2.2]{hovey2}).

We recall some standard definitions concerning
complexes of objects in a Gro\-then\-dieck category $\mathcal A$. Let $(M,\delta)$
(or just $M$, for simplicity) denote a chain complex in $\mathcal A$.
$$\cdots\to
M^{-1}\stackrel{\delta^{-1}}{\longrightarrow}
M^0\stackrel{\delta^{0}}{\longrightarrow}M^1\stackrel{\delta^{1}}{\longrightarrow}\cdots$$
We write $Z(M)=\cdots \to Z_nM\to Z_{n+1}M\to \cdots$ and
$B(M)=\cdots \to B_nM\to B_{n+1}M\to \cdots$ for the subcomplexes
consisting of the cycles and the boundaries of $M$.

Given an $M$ in $\mathcal A$, let $S^n(M)$
denote the complex which has $M$ in the $(-n)$th position and $0$ elsewhere ($n\in
\Z$). We denote by $D^n(M)$ the complex
$\cdots\to 0\to M\stackrel{id}{\to}M\to 0\to \cdots$ where $M$ is
in the $-(n+1)$th and $(-n)$th positions ($n\in \Z)$.

If $(M,\delta_M)$ and $(N,\delta_N)$ are two chain complexes, we
define $Hom(M,N)$ as the complex
$$\cdots\to \prod_{k\in \Z}\Hom(M^k,N^{k+n})\stackrel{\delta^n}
{\to}\prod_{k\in \Z}\Hom(M^k,N^{k+n+1})\to \cdots,$$ where
$(\delta^n f)^k={\delta}_N^{k+n}f^k-(-1)^nf^{k+1}{\delta}_M^k$.
Write ${\rm Ext}_{\Ch(\mathcal A)}(M,N)$ for the group of
equivalence classes of short exact sequences of complexes $0\to N\to
L\to M\to 0$. Let us note that $\Ch(\mathcal A)$ is a Grothendieck
category having the set $\{S^n(G):\ n\in \Z\}$ (or $\{D^n(G):\ n\in
\Z\}$) as a family of generators (where $G$ is a generator for
$\mathcal A$). So the functors ${\rm Ext}^i$, $i\in \Z$, can be
computed using injective resolutions.

Let $(\mathcal C, \Co)$ be a cotorsion pair in $\mathcal A$. We
recall from \cite{G0} the following definitions. An exact complex
$E$ in $\Ch(\mathcal A)$ is a \emph{$\Co$-complex} if $Z_nE\in \Co$,
for each $n\in \Z$. Let $\barCo$ denote the class of all
$\Co$-complexes. Then a complex $M=(M^n)$ in $\Ch(\mathcal A)$ is a
\emph{dg-$\mathcal C$} complex if $Hom(M,E)$ is an exact complex of
abelian groups for any complex $E\in \barCo$ and $M^n\in \mathcal
C$, for each $n\in \Z$. Let $\dgC$ denote the class of all
dg-$\mathcal C$ complexes of objects in $\mathcal A$.

Dually we can define the classes $\barC$ and $\dgCo$ of $\mathcal C$-complexes and
dg-$\Co$ complexes of objects in $\mathcal A$.

We will need the following lemma.

\begin{lemma}\label{exactk} Let $X$ be a scheme and $\kappa$ be a regular infinite
cardinal such that $\kappa \geq|V|$
and  $\kappa\geq|{\mathscr R}(v)|$ for all $v\in V$. Let $\mathscr N=(\mathscr
N^n),\mathscr M=(\mathscr M^n)$
be exact complexes of quasi--coherent sheaves on $X$ such that $\mathscr N\subseteq
\mathscr M$. For each $n \in \mathbb Z$, let ${\mathscr X}_n$ be a locally $\leq
\kappa$--presented quasi--coherent subsheaf of $\mathscr M ^n$. Then there exists an
exact complex of quasi--coherent sheaves $\mathscr T=(\mathscr T^n)$ such that
$\mathscr N\subseteq \mathscr T\subseteq M$, and for each $n \in \mathbb Z$,
$\mathscr T^n\supseteq \mathscr N^n+ \mathscr X_n$, and the quasi--coherent sheaf
$\mathscr T^n/\mathscr N^n$ is locally $\leq \kappa$--presented.
\end{lemma}
\noindent\Proof (I) First, consider the particular case of $\mathscr N = 0$. Let
${\mathscr Y}^n_{0} = {\mathscr X}_n + \delta^{n-1}({\mathscr X}_{n-1})$.
Then $({\mathscr Y}^n_0)$ is a subcomplex of $\mathscr M$.

If $i < \omega$ and ${\mathscr Y}^n_i$ is a locally $\leq \kappa$--presented quasi--coherent subsheaf of
$\mathscr M^n$, put ${\mathscr Y}^n_{i+1} = {\mathscr Y}^n_i + {\mathscr D}^n_i + \delta^{n-1}({\mathscr D}^{n-1}_i)$
where ${\mathscr D}^n_i$ is a locally $\leq \kappa$--presented quasi--coherent subsheaf of
$\mathscr M ^n$ such that $\delta^n ({\mathscr D}^n_i) \supseteq Z_{n+1}\mathscr M
\cap {\mathscr Y}^{n+1}_i$. (Such ${\mathscr D}^n_i$ exists by our assumption on
$\kappa$, since $Z_{n+1}\mathscr M \cap {\mathscr Y}^{n+1}_i \subseteq
\hbox{Ker}(\delta^{n+1}) = \mbox{Im}(\delta^n)$.)  Let ${\mathscr T}^n = \bigcup_{i
{<} \omega} {\mathscr Y}^n_i$. Then $Z_{n+1}\mathscr M \cap {\mathscr T}^{n+1} =
\bigcup_{i {<} \omega}(Z_{n+1}\mathscr M \cap {\mathscr Y}^{n+1}_i) \subseteq
\bigcup_{i {<} \omega} \delta^n ({\mathscr Y}^n_{i+1}) \subseteq  \delta^n ({\mathscr
T}^n)$. It follows that $\mathscr T=(\mathscr T^n)$ is an exact subcomplex of $\M$.
By our assumption on $\kappa$, $\mathscr T ^n$ is locally $\leq \kappa$--presented.

(II) In general, let $\bar{\mathscr M} = \mathscr M/\mathscr N$ and $\bar{\mathscr
X}_n = ({\mathscr X}_n + \mathscr N ^n)/\mathscr N ^n$. By part (I), there is
an exact complex of quasi--coherent sheaves $\bar{\mathscr T}$ such that
$\bar{\mathscr T} \subseteq \bar{\mathscr M}$, and
for each $n \in \mathbb Z$, $\bar{\mathscr T}^n\supseteq \bar{\mathscr X}_n$, and
the quasi--coherent sheaf $\bar{\mathscr T}^n$ is locally $\leq \kappa$--presented.
Then  $\bar{\mathscr T} = \mathscr T/\mathscr N$ for an exact subcomplex $\mathscr N
\subseteq \mathscr T \subseteq \mathscr M$, and
$\mathscr T$ clearly has the required properties.
\qed

\bigskip
As mentioned above, we will apply \cite[Theorem 2.2]{hovey2} to get a model
structure on
$\mathbb{C}(\Qco(X))$. We point out that $\Qco(X)$ is a closed symmetric monoidal
category under the tensor product (in the sense of \cite[Section 4.1]{Hov}) and
hence  $\mathbb{C}(\Qco(X))$ is also closed symmetric monoidal. We will therefore
investigate when the model structure is compatible with the induced closed symmetric
monoidal structure.

Let $X$ be a scheme with an associated quiver $Q_X$ (see Section \ref{filtqc}). Let
$\kappa \geq \mid V\mid$, and $\kappa \geq \mid \mathscr R(v)\mid$ for each $v\in V$.
We will assume that $X$ is \emph{semi--separated}, that is the
intersection of two affine open subsets of $X$ is again affine.

\medskip
\emph{For the rest of this section, we fix our notation as in Notation \ref{nttn}
and let $\lambda = \kappa^+$.}

We will moreover assume that $\mathcal C$ contains a generator of $\Qco(X)$. Then, by
Corollary \ref{und}, $(\mathcal C,{\mathcal L} ^{\perp})$ is a cotorsion pair.

\begin{lemma}\label{Ccomp}
$(\barC,{\mathit dg}\widetilde{{\mathcal L}^{\perp}})$ is a complete cotorsion pair
in $\Ch(\Qco(X))$.
\end{lemma}
\noindent\Proof $(\barC,{\mathit dg}\widetilde{{\mathcal
L}^{\perp}})$ is a cotorsion pair by  \cite[Corollary 3.8]{G0}.

We will prove that each complex $\mathscr C\in \barC$ is
$\widetilde{\mathcal L}$--filtered.  Then the completeness of
$(\barC,{\mathit dg}\widetilde{{\mathcal L}^{\perp}})$ follows as in
the proof of Corollary \ref{und} because $\barC$ contains a
generating set of $\Ch(\Qco(X))$ (for example $\{D^n(G) \mid n\in
\mathbb{Z}\}$ where $G\in \mathcal C$ is a generator of $\Qco(X)$).

Let $\mathscr C=(\mathscr M^n)\in \barC$. Then for each $n\in
\mathbb{Z}$, $Z_n\mathscr C\in \mathcal C$ and therefore
$Z_n\mathscr C$ has an $\mathcal L$--filtration $\mathcal
O_n=(\M^n_{\alpha}\mid \alpha\leq \sigma_n)$. For each
$n\in\mathbb{Z}$, $\alpha < \sigma_n$, consider a~locally
$\leq\kappa$--presented quasi--coherent sheaf $\mathscr A_{\alpha}^n
$ such that $\M_{\alpha+1}^n=\M_{\alpha}^n+\A_{\alpha}^n$, and the
corresponding family $\mathcal H_n$ as in Lemma \ref{HillQco}. Since
the complex $\mathscr C$ is exact, the $\mathcal L$--filtration
$\mathcal O_{n+1}$ determines a canonical prolongation of $\mathcal
O_n$ into a filtration $\mathcal O'_n=(\M_{\alpha}^n\mid \,
\alpha\leq \tau_n)$ of $\M^n$ where $\tau_n =  \sigma_n +
\sigma_{n+1}$ (the ordinal sum).

By definition, for each $\alpha \leq \sigma_{n+1}$, $\delta^n$ maps $\M_{\sigma_n + \alpha}^n$
onto $\M_\alpha^{n+1}$. So for each $\alpha < \sigma_{n+1}$ there is a
locally $\leq\kappa$--presented quasi--coherent subsheaf $\A^n_{\sigma_n + \alpha}$ of
$\M_{\sigma_n + \alpha+1}^n$ such that $\delta^n(\A^n_{\sigma_n + \alpha}) = \A^{n+1}_\alpha$.
Since for each $\sigma_n \leq \alpha < \tau_n$ we have $\hbox{Ker}(\delta^n) \subseteq \M_{\alpha}^n$,
it follows that $M_{\alpha+1}^n=\M_{\alpha}^n+\A_{\alpha}^n$.

Let $\mathcal H'_n$ be the family corresponding to $\A^n_\alpha$ ($\alpha < \tau_n$) by Lemma \ref{HillQco}. Since each closed subset of $\sigma_n$ is also closed when considered as a
subset of $\tau_n$, we have $\mathcal H_n\subseteq \mathcal H_n'$.
Note that ${Filt}(\mathcal L)\subseteq \mathcal C$, so $\mathcal
H_n'\subseteq \mathcal C$ by condition (3) of Lemma \ref{HillQco}.

Notice that $Z_n\mathscr C = \M ^n_{\sigma_n} = \sum_{\alpha < \sigma_n}\A_{\alpha}^n$.
We claim that for each closed subset $S\subseteq \tau_n$, we have
$Z_n\mathscr C \cap \sum_{\alpha\in S}\A_{\alpha}^n=
\sum_{\alpha \in S\cap \sigma_n}\A_{\alpha}^n\in \mathcal H_n$. To see this, we first
show that $\sum_{\alpha < \sigma_n}\A_{\alpha}^n(v)\cap \sum_{\alpha\in
S}\A_{\alpha}^n(v) = \sum_{\alpha \in S\cap \sigma_n}\A_{\alpha}^n(v)$ for each
$v\in V$. The inclusion $\supseteq$ is clear, so consider $a\in
(\sum_{\alpha < \sigma_n}\A_{\alpha}^n(v))\cap \sum_{\alpha\in S}\A_{\alpha}^n(v)$. Then
$a=a_{\alpha_0}+\cdots +a_{\alpha_k}$ where
$\alpha_i\in S$, $a_{\alpha_i}\in \mathscr A_{\alpha_i}^n(v)$ for all $i\leq k$, and
$\alpha_i > \alpha_{i+1}$ for all $i < k$. W.l.o.g., we can assume that $\alpha_0$ is
minimal possible. If $\alpha_0\geq\sigma_n$, then $a_{\alpha_0}=a-a_{\alpha_1}-
\cdots -a_{\alpha_k}\in (\sum_{\alpha < \alpha_0}\A_{\alpha}^n(v)) \cap \mathscr
A_{\alpha_0}^n(v)\subseteq \sum_{\alpha\in S,\alpha < \alpha_0} \mathscr
A_{\alpha}^n(v)$ as $\alpha_0\in S$ and $S$ is closed, in contradiction with the
minimality of $\alpha_0$. Hence $\alpha_0 < \sigma_n$, and $a\in \sum_{\alpha \in
S\cap \sigma_n}\A_{\alpha}^n(v)$. So $Z_n\mathscr C \cap \sum_{\alpha\in S}\A_{\alpha}^n =
\sum_{\alpha \in S\cap \sigma_n}\A_{\alpha}^n$, and the latter
quasi--coherent sheaf is in $\mathcal H_n$ because $S \cap \sigma_n$ is closed in
$\sigma_n$. This proves our claim.

By induction on $\alpha$, we will construct an $\widetilde{\mathcal L}$--filtration
$(\mathscr C_{\alpha}\mid \alpha\leq \sigma)$ of $\mathscr C$ such that $\mathscr
C_{\alpha}=(\mathscr N_{\alpha}^n)$, $Z_n\mathscr C_{\alpha}\in \mathcal H_n$ and
$\mathscr N_{\alpha}^n\in \mathcal H'_n$ for each $n \in \mathbb Z$.

First, $\mathscr C_0=0$, and if $\mathscr C_{\alpha}$ is defined and $\mathscr
C_{\alpha}\neq \mathscr C$, then for each $n\in \mathbb{Z}$ we take a locally $\leq
\kappa$--presented quasi--coherent sheaf $\mathscr X_n$ such that $\mathscr
X_n\nsubseteq \mathscr N_{\alpha}^n$ in case $\mathscr N_{\alpha}^n \subsetneq \M^n$
(this is possible by our assumption on $\kappa$), or $\mathscr X_n=0$ if $\M^n = \mathscr N_{\alpha}^n$.
If $\M^n = \mathscr N_{\alpha}^n$ for all $n \in \mathbb Z$, we let $\sigma = \alpha$
and finish our construction.

By Lemma \ref{exactk} there exists an exact subcomplex $\mathscr T=(\mathscr T^n)$
of $\mathscr C$ containing $\mathscr C_{\alpha}$ such that for each $n \in \mathbb
Z$, $\mathscr T^n\supseteq \mathscr N _\alpha^n+ \mathscr X_n$, and the
quasi--coherent sheaf $\mathscr T^n/\mathscr N _\alpha^n$ is locally $\leq
\kappa$--presented. Then $\mathscr Y_n=\mathscr T^n= \mathscr N_{\alpha}^n+\mathscr
X'_n$ for a locally $\leq \kappa$--presented quasi--coherent subsheaf $\mathscr
X'_n$ of $\mathscr M^n$. By condition (4) of Lemma \ref{HillQco} (for $\mathscr
N=\mathscr N^n_{\alpha}$ and $\mathscr X = \mathscr X_n'$), there exists a
quasi--coherent sheaf $\mathscr Y'_n=\mathscr P_n$ in $\mathcal H_n'$ such that
$\mathscr N_{\alpha}^n+\mathscr X_n'=\mathscr T^n\subseteq \mathscr P_n$ and
$\mathscr P_n/\mathscr N^n_{\alpha}$ is locally $\leq\kappa$--presented. Iterating
this process we obtain a countable chain $\mathscr Y_n\subseteq \mathscr
Y'_n\subseteq \mathscr Y''_n\subseteq \ldots$ whose union $\mathscr N^n_{\alpha +1} \in
\mathcal H'_n$ by condition (2) of Lemma \ref{HillQco}. Then $\mathscr C_{\alpha +1} =
(\mathscr N^n_{\alpha +1})$ is an exact subcomplex of $\mathscr C$ containing $\mathscr C_{\alpha}$.
Since $\mathscr N^n_{\alpha+1}\in \mathcal H'_n$, we have $Z_n\mathscr C_{\alpha+1}
= Z_n\mathscr C \cap \mathscr N^n_{\alpha+1} \in \mathcal H_n$ by the claim above.

In order to prove that $\mathscr C_{\alpha +1}/\mathscr C_{\alpha} \in \widetilde{\mathcal L}$, it remains to show
that for each $n \in \mathbb Z$, $Z_n(\mathscr C_{\alpha +1}/\mathscr C_{\alpha})\in \mathcal C$.
Since the complex $\mathscr C_{\alpha +1}/\mathscr C_{\alpha}$ is exact, it suffices to prove that
$F = (\delta^n (\mathscr N^n_{\alpha +1}) + \mathscr N^{n+1}_{\alpha})/\mathscr N^{n+1}_{\alpha} \in \mathcal C$.

We have $\mathscr N^n_{\alpha +1} = \sum_{\alpha\in S}\A_{\alpha}^n$ where w.l.o.g., $S$ is a closed subset of $\tau_n$ containing $\sigma_n$. Let $S^\prime = \{ \alpha < \sigma_{n+1} \mid \sigma_{n} + \alpha \in S \}$. Then $S^\prime$ is a closed subset on $\tau_{n+1} = \sigma_{n+1} + \sigma_{n+2}$. Indeed, for each $\alpha \in S^\prime$, we have
$$\sum_{\beta < \alpha}\A_{\beta}^{n+1}\cap \A_{\alpha}^{n+1} = \delta^n(\sum_{\beta < \sigma_n + \alpha} \A_{\beta}^{n}) \cap \delta^n(\A_{\sigma_n + \alpha}^{n}) \subseteq \delta^n(\sum_{\beta < \sigma_n + \alpha, \beta \in S} \A_{\beta}^{n}) = \sum_{\beta < \alpha, \beta \in S^\prime} \A_{\alpha}^{n+1}$$
where the inclusion $\subseteq$ holds because $S$ is closed in $\tau_{n}$ and $\hbox{Ker}(\delta ^n) \subseteq \sum_{\beta < \sigma_n + \alpha} \A_{\beta}^{n}$.

Since $\delta^n (\mathscr N^n_{\alpha +1}) = \sum_{\beta \in S^\prime} \A_{\beta}^{n+1}$, and $\mathscr N^{n+1}_{\alpha} = \sum_{\beta \in T} \A_{\beta}^{n+1}$ for a closed subset $T$ of $\tau_{n+1}$,
we have $F = \sum_{\beta \in S^\prime \cup T} \A_{\beta}^{n+1}/\sum_{\beta \in T} \A_{\beta}^{n+1}$, so $F \in \mathcal C$ by condition (3) of Lemma \ref{HillQco} for $\mathcal H^\prime_{n+1}$. This finishes the proof of $\mathscr C_{\alpha +1}/\mathscr C_{\alpha} \in \widetilde{\mathcal L}$.

If $\alpha$ is a limit ordinal we define $\mathscr C_{\alpha}=\bigcup_{\beta < \alpha}\mathscr C_{\beta}
= (\mathscr N_{\alpha}^n)$. Then $\mathscr N_{\alpha}^n\in \mathcal H'_n$ by condition (2) of
Lemma \ref{HillQco}, and $Z_n\mathscr C_{\alpha} = Z_n\mathscr C \cap \mathscr N^n_{\alpha}
\in \mathcal H_n$ by the claim above. This finishes the construction of the
$\widetilde{\mathcal L}$--filtration of $\mathscr C$.\qed

\medskip
Following \cite[Definition 6.4]{hovey2}, we call a cotorsion pair $(\mathcal
F,\mathcal C)$ in an abelian category $\mathcal A$ \emph{small} provided that (A1)
$\mathcal F$ contains a generator of $\mathcal A$, (A2) $\mathcal C = \mathcal S
^\perp$ for a subset $\mathcal S \subseteq \mathcal F$, and (A3) for each $S \in
\mathcal S$ there is a monomorphism $i_S$ with cokernel $S$ such that if $\mathcal
A(i_S,X)$ is surjective for all $S \in \mathcal S$, then $X \in \mathcal C$.

We now show that condition (A3) above is redundant in case $\mathcal A$ is a
Grothendieck category:

\begin{lemma}\label{small2}
Let $(\mathcal F,\mathcal C)$ be a cotorsion pair in a Grothendieck category
$\mathcal A$ satisfying conditions (A1) and (A2) above. Then $(\mathcal F,\mathcal
C)$ is  small.
\end{lemma}
\noindent\Proof We will show that $(\mathcal F,\mathcal C)$ satisfies a slightly
weaker version
of condition (A3), namely that for each $L\in {\mathcal S}$ there is a set $\mathcal
E _L$ of
exact sequences $0\to K\to U\to L\to 0$ such that $Y\in \mathcal C$ if and only if
$\Hom(U,Y)\to
\Hom(K,Y)\to 0$ is exact for each exact sequence in $\mathcal E _L$. For a given
$L$, we define $\mathcal E _L$ as the set
of all representatives of short exact sequences $0\to K\to U\to L\to 0$ where $U$ is
$\leq \kappa$--presented (where
$\kappa$ comes from \cite[Corollary 2.3]{EEGR0} for $Y = L$; in particular, we can
take $\kappa$ is as in Notation \ref{nttn}
in case $\mathcal A = \Ch(\Qco(X))$).

Suppose that $G$ is an object of $\mathcal A$ such that $\Hom(U,G)\rightarrow
\Hom(K,G) \rightarrow 0$ is exact for each
exact sequence in $\mathcal E _L$. We will prove that ${\rm Ext}^1(L,G)=0$ for all
$L \in \mathcal F$. By condition (A2),
it suffices to prove that ${\rm Ext}^1(L,G)=0$ for all $L\in {\mathcal S}$. So let
$0\rightarrow G\rightarrow V\rightarrow L\rightarrow 0$ be exact
with $L\in {\mathcal S}$. We want to show that this sequence splits.
By our choice of $\kappa$, there is $U\subseteq V$ such that $U$ is $\leq
\kappa$--presented and $V = G + U$.
Then the sequence $0\rightarrow G\cap U \rightarrow U\rightarrow L \rightarrow 0$ is
isomorphic to one in $\mathcal E _L$.

Consider the commutative diagram

$$\begin{CD}
0@>>> {G\cap U}@>>> U@>>> L @>>> 0\\
@. @VVV @VVV @| @.\\
0@>>> {G}@>>> {G+U} @>>> L @>>> 0.\\
\end{CD}$$

Our hypothesis now implies that the inclusion $G\cap U \rightarrow G$ can be
extended to $U\rightarrow G$ so,
since the left--hand square is a pushout, we see that the bottom row splits.
This proves that ${\rm Ext}^1(L,G)=0$. Now, replacing the set $\mathcal S$ by
$\mathcal S ^\prime =
\{ L^{(card(\mathcal E _L))} \mid L \in \mathcal S \}$, we see that both conditions
(A2) and (A3)  hold for $\mathcal S ^\prime$,
hence  the cotorsion pair $(\mathcal F,\mathcal C)$ is small.
\qed

\medskip
Now we can prove the main theorem of our paper.

\begin{theorem}\label{model}
Let $X$ be a semi-separated scheme. There is a model category
structure on $\mathbb{C}(\Qco(X))$ where the weak equivalences are
the homology isomorphisms, the cofibrations (resp.\ trivial
cofibrations) are the monomorphisms with cokernels in $\dgC$ (resp.\
in $\barC$), and the fibrations (resp.\ trivial fibrations) are the
epimorphisms whose kernels are in $\dgCoo$ (resp.\ $\barCoo$).

Moreover, if every
$M\in S_v$ is a flat $\mathscr R(v)$-module, and $M\otimes_{\mathscr R(v)} N\in S_v$
whenever $M, N\in S_v$, then the model structure is monoidal with respect to the
usual tensor product of complexes of quasi-coherent sheaves.
\end{theorem}
\noindent\Proof We will apply Hovey's Theorem \cite[Theorem 2.2]{hovey2}. First, the
results of \cite[Section 5]{hovey2} guarantee that the weak equivalences of our
model structure are the homology isomorphisms. In our case $\mathcal W$ is the class
of all exact complexes of quasi--coherent sheaves. It is easy to check that this is
a thick subcategory of $\Ch(\Qco(X))$. Now, according to Hovey's Theorem, we will
have to check that the pairs $(\dgC,\dgCoo\cap \mathcal W)$ and $(\dgC\cap \mathcal
W,\dgCoo)$ are complete cotorsion pairs (notice that our notion of completeness
coincides with Hovey's notion of 'functorial completeness'). We will proceed in
three steps, proving that
\begin{enumerate}
\item The pairs $(\barC,\dgCoo)$ and $(\dgC,\barCoo)$ are cotorsion pairs.
\item $\dgC\cap \mathcal W=\barC$ and $\dgCoo\cap \mathcal W=\barCoo$.
\item The cotorsion pairs $(\barC,\dgCoo)$ and $(\dgC,\barCoo)$ are complete.
\end{enumerate}

Condition $(1)$ follows from \cite[Corollary 3.8]{G0}.

\medskip
Let us check condition $(2)$. By \cite[Corollary 3.9]{G} (4. $\Rightarrow$ 1.) it
suffices to prove that $\dgC\cap \mathcal W=\barC$. The inclusion $\barC\subseteq
\dgC\cap\mathcal W$ was proven in \cite[Lemma
3.10]{G0}.
Let us prove that $\dgC\cap \mathcal W\subseteq \barC$.
So let ${\mathscr Y}$ be a complex in $\dgC\cap \mathcal W$ (so ${\mathscr
Y}(v)$ is a complex of $\R(v)$-modules, for all $v\in V$).
To see that $\mathscr Y$ is in $\barC$
we have to check that $Z_n{\mathscr Y}\in \mathcal C$, for all $n\in \mathbb{Z}$.
But this means that the $\R(v)$-module $Z_n{\mathscr Y}(v)$ belongs
to ${\mathcal F}_v$ for all $v\in V$. By \cite[Corollary 3.9]{G} if a complex of
$\R(v)$-modules is exact and belongs to $\dgFv$ then it belongs to $\barFv$ (so
$Z_n{\mathscr Y}(v)\in {\mathcal F}_v$ for all $v\in V$). Therefore we will be done
if we prove that $\mathscr Y(v)$ is exact and belongs to $\dgFv$. Since the complex
$\mathscr Y$ is exact, for each affine open set $v\in V$, $\mathscr Y(v)$ is an
exact complex of
$\R(v)$-modules. Let us see that $\mathscr Y(v)\in \dgFv$, for all $v\in V$. So let
$E$ be a complex of $\R(v)$-modules in $\barSo$ (so $E$ is exact and $Z_nE\in
{\mathcal S_v}^{\perp}$). We have to check that $ Hom(\mathscr Y(v),E)$
is exact. Since $X$ is semi-separated, by \cite[Proposition 5.8]{HAR} there exists a
right adjoint $i_*{_v}:\R(v)\Mod\to \Qco(X)$ of the restriction functor
$i^*{_v}:\Qco(X)\to \R(v)\Mod$ (defined by $i^*{_v}(\mathscr M)=\mathscr M(v)$). The
adjointness situation can be lifted up to $\Ch(\Qco(X))$. Then there is an
isomorphism $$Hom_{\Ch(\R(v))}(\mathscr Y(v),E)=Hom_{\Ch(\R(v))}(i^*{_v}(\mathscr
M),E)\cong Hom_{\Ch(\Qco(X))}(\mathscr Y,i_*{_v}(E)) $$ and since the functor
$i_*{_v}$ preserves exactness, $i_*{_v}(E)$ will be an exact complex in
$\Ch(\Qco(X))$. Since $\mathscr Y\in \dgC$, once we show that $i_*{_v}(E)\in \barLo$
we will finish by the comment above. But, $Z_ni_*{_v}(E)=i_*{_v}(Z_nE)$. Hence, for
each $\mathscr T\in \mathcal C$,
$${\rm Ext}^1_{\Qco(X)}(\mathscr T, i_*{_v}(Z_nE))\cong {\rm
Ext}^1_{\R(v)}(i^*{_v}(\mathscr T), Z_nE)=0,$$where the last equality
follows because $i^*{_v}(\mathscr T)=\mathscr T(v)\in \mathcal F_v$ and $Z_nE\in
{\mathcal S_v}^{\perp}$.

\medskip
Now let us prove condition $(3)$. By Lemma \ref{Ccomp} the cotorsion pair
$(\barC,\dgCoo)$ is complete. We claim that the cotorsion pair $(\dgC,\barCoo)$ is
also complete. Let $\mathcal I ^\prime$ be a set of representatives of the
quasi--coherent sheaves in $\mathcal L$. Then clearly
$(\mathcal I ^\prime)^\perp =\Uo$. We will prove that $\mathcal{I}^\perp = \barCoo$
where $\mathcal{I}=\{S^n(\A)\, \mid\, \A\in \mathcal I^\prime,\,n\in \Z\}\cup
\{S^n(\mathscr G)\,\mid\, n\in \Z\}$ (and $\mathscr G\in \mathcal C$ is a generator
of $\Qco(X)$). Then the claim will follow by Lemma \ref{small2} and \cite[Corollary
6.6]{hovey2}. It is easy to check that
$\mathcal{I}\subseteq \dgC$ for $S^m(\mathscr A)^l\in \mathcal C$ ($l\in \Z$),
and for every exact complex $\mathscr M\in \barCoo$, $Hom(S^m(\A),\mathscr M)$ is
the complex
$$\cdots\to \Hom(\A,\mathscr M^l)\to
\Hom(\A,\mathscr M^{l+1})\to\cdots$$which is obviously exact because
$Z_n\mathscr M, B_n\mathscr M\in \Uo$. Therefore $\mathcal{I}^{\perp}\supseteq
(\dgC)^{\perp}=\barCoo$. We now prove the converse: let $\mathscr N\in
\mathcal{I}^{\perp}$. We have to see that $\mathscr N$ is exact and that
$Z_n\mathscr N\in \Uo$. First, we prove that $\mathscr N$ is exact. It is clear that
this is
equivalent to each morphism $S^n(\mathscr G)\to \mathscr N$ (for $\mathscr G$ a
generator of $\Qco(X)$) being
extendable to $D^n(\mathscr G)\to \mathscr N$ for each $n\in \Z$. But
this follows from the short exact sequence
$$0\to S^n(\mathscr G)\to D^n(\mathscr G) \to S^{n+1}(\mathscr G) \to 0 $$
since ${\rm Ext}^1(S^{n+1}(\mathscr G),\mathscr N)=0$. Now we prove that
$Z_n\mathscr N\in \Uo$.
Since $\mathcal I ^\perp = \Uo$ we only need to prove that ${\rm
Ext}^1_{\Qco(X)}(\A,Z_n \mathscr N)=0$ for all $n\in \Z$. But there exists a
monomorphism of abelian groups $$0\to {\rm Ext}^1_{\Qco(X)}(\A,Z_n \mathscr N)\to
{\rm Ext}^1_{\Ch(\Qco(X))}(S^{-n}(\A),\mathscr N)$$  (see e.g.\ \cite[Lemma
5.1]{EEGR}) and since the latter group is 0, we get that $Z_n \mathscr N\in \Uo$.
This proves our claim, and thus finishes the proof of condition $(3)$.

\bigskip
Finally to get that the model structure is monoidal we apply \cite[Theorem 5.1]{G}
(by noticing that the argument of the proof of \cite[Theorem 5.1]{G} carries over
without the assumption of $\mathcal F$ being closed under direct limits).
If $S_v$ is contained in the class of all flat modules then every quasi--coherent
sheaf in $\mathcal C$ is flat. So condition (1) of \cite[Theorem 5.1]{G} holds. Now
if $M\otimes_{\mathscr R(v)}N\in S_v$, where $M,N$ are $\mathscr R(v)$-modules in
$S_v$, it follows that $L\otimes_{\mathscr R(v)}T \in \mathcal F_v$, where $L,T$ are
$S_v$-filtered $\mathscr R(v)$-modules (because the tensor product commutes with
direct limits). And so $\mathscr L\otimes_{\mathscr R}\mathscr T\in \mathcal C$, for
any ${\mathscr L},{\mathscr T}\in \mathcal C$. So condition $(2)$ of \cite[Theorem
5.1]{G} also holds. Finally condition (3) of \cite[Theorem 5.1]{G} is immediate
because, for all $v\in V$, $\mathcal F_v$ contains all projective $\mathscr
R(v)$-modules, so in particular $\mathscr R\in \mathcal C$.\qed

\bigskip
\emph{The proof of Corollaries \ref{cor1} and \ref{cor2}.} In Theorem \ref{model}, we
take $S_v = \{ R(v) \}$, and $S_v$ = a representative set of all flat modules of
cardinality
$\leq \card (R(v)) + \aleph_0$, respectively. Notice that in the first case,
$\mathcal C$
is the class all of vector bundles, while in the second, $\mathcal C$ is the class
of all flat
quasi--coherent sheaves. \qed

\medskip
If $X={\mathbf P}^n(R)$ where $R$ is any commutative noetherian ring, then every
quasi--coherent sheaf on $X$ is a filtered union of coherent subsheaves, and the
family of so--called twisting sheaves $\{\mathscr O(n)\mid n\in \mathbb{Z}\}$
generates the category of coherent sheaves on $X$ cf.\ \cite[Corollary 5.18]{HAR},
so $\bigoplus_{i\in \mathbb{Z}}\mathscr O(n)$ is a (vector bundle) generator for
$\Qco(X)$.  So Corollary \ref{cor1} applies to this setting. In particular, we
extend here \cite[Theorem 6.1]{EEGR} which deals with the case of the projective
line.

\medskip
Finally, we consider the case of restricted Drinfeld vector bundles:

\medskip
\emph{The proof of Corollary \ref{cor3}.} In view of Theorem \ref{model},
the proof will be complete once we show

\begin{lemma}
If $R$ is a commutative ring and $M$ and $N$ are $\leq \kappa$-presented flat
Mittag-Leffler modules, then so is $N\otimes_R M$.
\end{lemma}
\noindent\Proof It is clear that $N\otimes_R M$ is $\leq \kappa$-presented. Let us
check the Mittag-Leffler condition (see Definition \ref{flml}). So let $( M_i \mid i
\in I )$ be a family of $R$-modules. Since $N$ is flat Mittag-Leffler the canonical
map
$N \otimes_R \prod_{i \in I} M_i \to \prod_{i \in I} N \otimes_R M_i$ is
a monomorphism. Now since $M$ is flat, we get a monomorphism $$(M\otimes_R N)
\otimes_R \prod_{i \in I} M_i\cong M\otimes_R (N \otimes_R \prod_{i \in I} M_i) \to
M\otimes_R (\prod_{i \in I} N \otimes_R M_i).$$ Now we apply the fact that $M$ is
Mittag-Leffler to the family $( N\otimes_R M_i \mid i \in I )$ to get a monomorphism
$$M\otimes_R (\prod_{i \in I} N \otimes_R M_i)\to \prod_{i\in I} M\otimes_R
(N\otimes_R M_i)\cong \prod_{i\in I} (M\otimes_R N)\otimes_R M_i. $$So the claim
follows by composing the previous monomorphisms. \qed

\section{Flat Mittag--Leffler Abelian Groups}

Let $X$ be a scheme having a generating set consisting of Drinfeld
vector bundles. We have already seen that restricted Drinfeld vector
bundles impose monoidal model structures on $\mathbb{C}(\Qco(X))$
whose weak equivalences are the homology isomorphisms (see Corollary
\ref{cor3}). This result suggests that the entire class of all
Drinfeld vector bundles could also impose a cofibrantly generated
model structure in $\mathbb{C}(\Qco(X))$. The aim of this section is
to prove Theorem \ref{notmodel} and show thus that this is not the
case in general.

We recall that, given a class of objects $\mathcal F$ in a Grothendieck category
$\mathcal A$, an \emph{$\mathcal F$-precover} of an object $M$
is a morphism $\varphi:F\to M$ with $F\in \mathcal F$ such that
${\rm Hom}_{\mathcal A}(F',F)\to {\rm Hom}_{\mathcal A}(F',M)$ is an epimorphism for
every $F'\in\mathcal F$.
The class $\mathcal F$ is said to be precovering if every object of $\mathcal A$
admits an $\mathcal F$-precover (see \cite[Chapters 5 and 6]{EdO} for properties of
such classes).
For example the class of projective modules $\mathcal P$ is precovering.
Similarly as $\mathcal P$ is used to define projective resolutions, one can employ a
precovering class $\mathcal F$ to define $\mathcal F$--resolutions and a version of
relative homological algebra can be developed (see \cite{EdO}).

\medskip
Let $\mathcal D$ denote the class of all flat Mittag--Leffler modules over a ring $R$.

Both the flat and the Mittag--Leffler modules are clearly closed under
pure--submodules,
hence so is the class $\mathcal D$. Similarly, $\mathcal D$
is closed under (pure) extensions (see \cite[2.1.6]{RG}).
Moreover, the countably generated modules in $\mathcal D$
are exactly the countably generated projective modules by \cite[2.2.2]{RG}.
In general, the modules in $\mathcal D$ are characterized
by the following theorem from \cite{RG} (cf.\ \cite{D}):

\begin{theorem}\label{ML}
Let $R$ be a ring and $M \in \rmod R$. Then the following are equivalent:
\begin{enumerate}
\item $M$ is a flat Mittag-Leffler module.
\item Every finite (or countable) subset of $M$ is contained in
a countably generated projective submodule which is pure in $M$.
\item Every finite subset of $M$ is contained in a projective submodule which is
pure in $M$.
\end{enumerate}
\end{theorem}

This theorem implies that the class $\mathcal D$ is closed under $\mathcal
D$--filtrations
(see Definition \ref{defilt}).

We start with a more specific characterization of flat Mittag--Leffler modules in
the particular
case of Dedekind domains with countable spectrum.
Recall that a module $M$ over a right hereditary ring is {\em
$\aleph_1$--projective} ({\em $\aleph_1$--free}) provided
that each countably generated submodule of $M$ is projective (free).

\begin{lemma}\label{dedek}
Assume that $R$ is a Dedekind domain such that $\mbox{Spec}(R)$ is countable.
Let $M \in \rmod R$. Then $M \in \mathcal D$ iff $M$ is $\aleph_1$--projective.
\end{lemma}
\noindent\Proof If $M \in \mathcal D$ and $C$ is a countably
generated submodule of $M$, then $C$ is contained in a countably generated projective
(and pure) submodule $P$ of $M$ by Theorem \ref{ML}. Since $R$ is right hereditary,
$C$ is also projective.

In order to prove the converse, in view of Theorem \ref{ML},
it suffices to prove that each countable subset $C$ is contained
in a countably generated pure submodule $P$ of $M$. Since $M$ is
flat and $P$ is projective, the purity of the embedding $P \subseteq M$
can be tested only w.r.t.\ all simple modules by \cite[Lemma 11]{GH}, that is,
one only has to construct a countably generated module $P \supseteq C$ such
that $P \otimes_R S \to M \otimes_R S$ is monic for each simple left module $S$.
But the class of all simple modules has a countable set of representatives by
assumption,
so we obtain our claim by applying \cite[I.8.8]{S}.
\qed

\medskip
From now on, we will restrict ourselves to the particular case of (abelian) groups.
By Lemma \ref{dedek}, a group $A$ is flat and Mittag--Leffler iff $A$ is
$\aleph_1$--free.
Our aim is to show that the class of all $\aleph_1$--free
groups is not precovering. We will prove this following an idea from \cite{ES}
where the analogous result was proven consistent with (but independent of) ZFC + GCH
for the subclass of $\mathcal D$ consisting of all Whitehead groups.

The reason why our result on $\mathcal D$ holds in ZFC
rather than only in some of its forcing extensions rests in
the following fact whose proof goes back to \cite{H} (see also \cite{GT0}):
for each non--cotorsion group $A$, there is a  {\em Baer--Specker group} (that is,
the product $\mathbb Z
^\kappa$ for some $\kappa \geq \aleph_0$) such that $\Ext 1{\mathbb Z}{\mathbb Z
^\kappa}A
\neq 0$;
moreover, the Baer--Specker group can be taken small in the following sense:

\begin{lemma}\label{cotors} Define a sequence of cardinals
$\kappa_\alpha$ ($\alpha \geq 0$) as follows:
\begin{itemize}
\item $\kappa_0 = \aleph_0$,
\item $\kappa_{\alpha + 1} = \sup_{i {<} \omega} \kappa_{\alpha,i}$
where $\kappa_{\alpha,0} = \kappa _\alpha$ and $\kappa_{\alpha,n+1} =
2^{\kappa_{\alpha,n}}$, and
\item $\kappa_\alpha = \sup_{\beta {<} \alpha} \kappa_\beta$ when $\alpha$ is a limit
ordinal.
\end{itemize}
Let $\alpha$ be an ordinal and $A$ be a non--cotorsion group of cardinality $\leq
2^{\kappa_\alpha}$.
Then $\Ext 1{\mathbb Z}{\mathbb Z ^{\kappa_\alpha}}A \neq 0$.
\end{lemma}
\noindent\Proof This is a consequence of \cite[1.2(4)]{GT0} where the following
stronger assertion is proven:

'If $A$ is any group of cardinality $\leq 2^{\kappa_\alpha}$
such that $\Ext 1{\mathbb Z}{D_{\kappa_\alpha}}A = 0$
(where $D_{\kappa_\alpha}$ is a certain subgroup of $\mathbb Z ^{\kappa_\alpha}$)
then $\Ext 1R{\mathbb Q}A = 0$, that is, $A$ is a cotorsion group.'
\qed

\begin{remark}
Under GCH, the definition of the $\kappa_{\alpha}$'s simplifies as follows: if
$\kappa_{\alpha}=\aleph_{\beta}$, then $\kappa_{\alpha+1}=\aleph_{\beta+\omega}$.
\end{remark}

It follows that though the class $\mathcal D$ of all $\aleph_1$--free groups is
closed under $\mathcal D$--filtrations, it is not of the form $^\perp \mathcal C$
for any class of groups $\mathcal C$:

\begin{theorem}\label{corol1} Let $R = \mathbb Z$. Then $\mathcal D \neq {}^\perp
\mathcal C$, for each class
$\mathcal C \subseteq \rmod{\mathbb Z}$. In fact, $^\perp (\mathcal D ^\perp)$ is
the class of all flat (= torsion-free) groups.
\end{theorem}
\noindent\Proof Since all the Baer--Specker groups are $\aleph_1$--free (see e.g.\
\cite[IV.2.8]{EM}),
Lemma \ref{cotors} implies that $\mathcal D ^\perp$ coincides with the class of all
cotorsion groups, so
$^\perp (\mathcal D ^\perp)$ is the class of all flat groups. Since $\mathbb Q
\notin \mathcal D$,
we have $\mathcal D \neq {}^\perp \mathcal C$ for each class $\mathcal C \subseteq
\rmod{\mathbb Z}$.
\qed

\begin{lemma} \label{corol2} Let $\alpha$ be an ordinal and $A$ be an
$\aleph_1$--free group
of cardinality $\leq 2^{\kappa_\alpha}$ where $\kappa_\alpha$ is defined as in Lemma
\ref{cotors}.
Then $\Ext 1{\mathbb Z}{\mathbb Z ^{\kappa_\alpha}}A \neq 0$.
\end{lemma}
\noindent\Proof  In view of Lemma \ref{cotors},
it suffices to verify that no non--zero $\aleph_1$--free group is cotorsion.
Indeed, each reduced torsion free cotorsion group $A$ has a direct
summand isomorphic to $\mathbb J _p$ (the group of all $p$--adic
integers for some prime $p \in \mathbb Z$) by \cite[V.2.7 and V.2.9(5),(6)]{EM}.
However, if $A$ is $\aleph_1$--free, then it is cotorsion--free by
\cite[V.2.10(ii)]{EM},
so $\mathbb J _p$ does not embed into $A$.
\qed

\begin{remark}\label{nofilt} Theorem \ref{corol1} already implies that we cannot
improve
Corollary \ref{cor3} by extending the claim to all Drinfeld vector bundles (that is,
removing the $\kappa$-filtration restriction): consider the afine scheme $X = \Spec
(\mathbb Z)$. Then there is a category equivalence $\Qco (X) \cong \rmod {\mathbb
Z}$ by \cite[Corollary 5.5]{HAR}. By Theorem \ref{corol1}, for each infinite
cardinal $\kappa$, there is a Drinfeld vector bundle $\mathscr M$ which does not
have a $\mathcal C$-filtration where $\mathcal C$ is the class of all locally $\leq
\kappa$-presented Drinfeld vector bundles.

In more detail, if $\kappa=\kappa_{\alpha}$ (see Lemma \ref{cotors}) and $\mathcal D
^{\leq \kappa}$ denotes the class of all $\leq \kappa$-generated
$\aleph_1$-free groups, then $\mathbb Z ^{2^\kappa} \in \mathcal D
\setminus {}^\perp((\mathcal D ^{\leq \kappa})^\perp)$.

Indeed, denote by $\mathcal E$ a representative set of elements of the
class $\mathcal D ^{\leq \kappa}$. Then $\mid \mathcal E \mid \leq
2^\kappa$, so by \cite[Theorem 2]{ET}, there exists $A \in \mathcal E
^\perp$ such that $\mid \mathcal A \mid = 2^{2^\kappa}$ and $A$ has a
$\mathcal E$-filtration. In particular, $A$ is $\aleph_1$-free, and $\Ext
1{\mathbb Z}{\mathbb Z ^{2^\kappa}}A \neq 0$ by Lemma \ref{corol2}. Hence $\mathbb
Z ^{2^\kappa} \notin  {}^\perp(\mathcal E ^\perp)$.
\end{remark}

In view of Remark \ref{nofilt}, the class of all $\aleph_1$--free groups cannot induce
a cofibrantly generated model category structure on
$\Qco(\rm{Spec}(\mathbb{Z}))\cong \rm{Mod-}\mathbb{Z}$
compatible with its abelian structure. This is the second claim of Theorem
\ref{notmodel}.
In order to prove the (stronger) first claim of Theorem \ref{notmodel}, it remains
to show the following

\begin{theorem}\label{noprec} The class of all $\aleph_1$--free groups is not
precovering.
\end{theorem}

\noindent\Proof Assume there exists a $\mathcal D$--precover of $\mathbb Q$,
and denote it by $f : B \to \mathbb Q$. We will construct an $\aleph_1$--free group
$G$
of infinite rank such that there is no non--zero homomorphism from $G$ to $B$. Since
$G$
 has infinite rank and $\mathbb Q$ is injective, there is a non--zero (even
surjective)
 homomorphism $g : G \to \mathbb Q$. Clearly $g$ does not factorize through $f$, a
contradiction.

First, we take an ordinal $\alpha$ such that $\mu = 2^{\kappa_\alpha} \geq \card (B)$
(see Lemma \ref{cotors}). The $\aleph_1$--free group $G$ will be the last
term of a continuous chain of $\aleph_1$--free groups of infinite rank,
$( G_\nu \mid \nu \leq \tau )$, of length $\tau \leq \mu^+$.
The chain will be constructed by induction on $\nu$ as follows: first, $G_0$ is any
free group of infinite rank.

Assume $G_\nu$ is defined for some $\nu {<} \mu^+$ and consider the set $I_\nu$
of all non--zero homomorphisms from $G_\nu$ to $B$. If $I_\nu = \emptyset$,
we put $\tau = \nu$ and finish the construction.
Otherwise, we fix a free presentation $0 \to K \hookrightarrow F \to \mathbb Z
^{\kappa_\alpha} \to 0$
of $\mathbb Z ^{\kappa_\alpha}$, and denote by $\theta$ the inclusion of $K$ into $F$.

For each $h \in I_\nu$, let $A_h$ be the image of $h$. By Lemma \ref{corol2},
$\Ext 1{\mathbb Z}{\mathbb Z ^{\kappa_\alpha}}{A_h} \neq 0$, so there exists a
homomorphism $\phi_h : K \to A_h$ which does not extend to $F$. Since $K$ is
free and $h$ maps onto $A_h$, there is a homomorphism $\psi_h : K \to G_\nu$ such
that $h \psi_h = \phi_h$.

Denote by $\Theta$ the inclusion of $K^{(I_\nu)}$ into $F^{(I_\nu)}$, and define
$\Psi \in \OHom {\mathbb Z}{K^{(I_\nu)}}{G_\nu}$ so that the $h$-th component of
$\Psi$
is $\psi_h$, for each $h \in I_\nu$.

The group $G_{\nu + 1}$ is defined by the pushout of $\Theta$ and $\Psi$:

$$\begin{CD}
K^{(I_\nu)} @>{\Theta}>>  F^{(I_\nu)} \\
@V{\Psi}VV  @V{\Omega}VV \\
G_\nu @>{\subseteq}>> G_{\nu + 1}.
\end{CD}$$

Note that $G_{\nu + 1}/G_\nu \cong F^{(I_\nu)}/K^{(I_\nu)}$
 is $\aleph_1$--free because $\mathbb Z ^{\kappa_\alpha}$
 is $\aleph_1$--free by \cite[IV.2.8]{EM}. It follows that $G_{\nu +1}$
 is an $\aleph_1$--free group of infinite rank.

If $\nu \leq \mu^+$ is a limit ordinal we put $G_\nu = \bigcup_{\sigma {<} \nu}
G_\sigma$.
 Clearly $G_\nu$ has infinite rank, and since $G_{\sigma + 1}/G_\sigma$ is
$\aleph_1$--free
 for each $\sigma {<} \nu$ by construction, $G_\nu$ is also $\aleph_1$--free.

It remains to show that there exists $\nu \leq \mu^+$ such that $I_\nu = \emptyset$.
Assume $I_\nu \neq \emptyset$ for all $\nu {<} \mu^+$ (hence $G_\nu$ is defined for all
 $\nu \leq \mu^+$); we will prove that $I_{\mu^+} = \emptyset$.

Assume there is a non--zero homomorphism $f : G_{\mu^+} \to B$ and let $\nu {<} \mu^+$
 be such that $h := f \restriction G_\nu \neq 0$.

Using the notation introduced in the non--limit step of the construction,
we will prove that $A_h$ is a proper submodule of the image of
$h^\prime = f \restriction G_{\nu + 1}$. If not, then $h^\prime \Omega$ extends $h
\Psi$
to a homomorphism $F^{(I_\nu)} \to A_h$. Denote by $\iota_h$ and $\iota^\prime_h$
the $h$-th
canonical embedding of $K$ into $K^{(I_\nu)}$ and of $F$ into $F^{(I_\nu)}$,
respectively.
Then $h^\prime \Omega \iota^\prime_h$ extends $h \Psi \iota_h = h \psi_h = \phi_h$
to a
homomorphism $F \to A_h$, in contradiction with the definition of $\phi_h$.

This proves that the image of $f \restriction G_\nu$ is a proper submodule of the
image of
$f \restriction G_{\nu + 1}$ for each $\nu \in C$, where $C$ is the set of all $\nu
{<} \mu^+$
such that $f \restriction G_\nu \neq 0$. However, $f \neq 0$ implies that $C$ has
cardinality $\mu^+$,
in contradiction with $\card (B) {<} \mu^+$.
This proves that $\OHom {\mathbb Z}{G_{\mu^+}}B = 0$, that is, $I_{\mu^+} =
\emptyset$.
\qed

\end{document}